\documentclass{article}

\usepackage{graphicx}
\usepackage{amsmath,amsfonts}

\newtheorem{theorem}{Theorem}
\newtheorem{acknowledgement}[theorem]{Acknowledgement}

\newtheorem{definition}[theorem]{Definition}
\newtheorem{example}[theorem]{Example}

\newtheorem{lemma}[theorem]{Lemma}

\newtheorem{proposition}[theorem]{Proposition}
\newtheorem{remark}[theorem]{Remark}

\newenvironment{proof}[1][Proof]{{\sc #1.} }{\hfill $\Box$}

\begin{document}

\title{On a relation between stochastic integration and geometric measure theory}
\author{
\\
  { Franco Flandoli and Massimiliano Gubinelli\footnote{Corresponding
      author, e-mail {\tt m.gubinelli@dma.unipi.it}} }              \\
  {\small\it Dipartimento di Matematica Applicata,
Universit\`{a} di Pisa}    \\[-0.1cm]
  {\small\it Via Bonanno 25B, 56126 Pisa, Italia}          \\[-0.1cm]
  \\[-0.1cm]  \and
  { Mariano Giaquinta }              \\
  {\small\it Scuola Normale Superiore di Pisa }    \\[-0.1cm]
  {\small\it P.za dei Cavalieri  7, 56125 Pisa, Italia}          \\[-0.1cm]
  \\[-0.1cm]  \and
  { Vincenzo M. Tortorelli }              \\
  {\small\it  Dipartimento di Matematica,
Universit\`{a} di Pisa }    \\[-0.1cm]
  {\small\it Via Buonarroti  2, 56127 Pisa, Italia}          \\[-0.1cm]
}
\date{November 2002}
\maketitle

\begin{abstract}
Two problems are addressed for the path of certain stochastic processes: a)
do they define currents? b) are these currents of a classical type? A
general answer to question a) is given for processes like semimartingales or
with Lyons-Zheng structure. As to question b), it is shown that H\"{o}lder
continuous paths with exponent $\gamma > 1/2$
define integral flat chains.
\end{abstract}

\section{Introduction}

This work starts from two motivations. The first one is the intuition
by one of the authors 
 and T. Lyons that interesting relations could exist between the
theory of currents and the theory of rough paths (see \cite{GMS} and \cite
{Lyons} as basic references on these topics). The second motivation comes
from recent investigations on stochastic models of vortex
filaments in turbulent fluids, where the vorticity field is rather singular
and a description by means of current looks promising (see \cite{Ffil}, \cite
{FGub}, \cite{Ascona} and the previous works \cite{Ch}, \cite{Ga}, 
\cite{PLLMaj}). The present work provides only some initial results on these
two problems.

We address the following two questions:

\vspace{0.5cm}

a)\ do the paths of some stochastic processes define $1$-currents?

\vspace{0.5cm}

b) are these currents of a classical type?

\vspace{0.5cm}

The definition of random current is given in the next section, but roughly speaking
our subject of investigation are the Stratonovich stochastic integrals of
the form 
\begin{equation*}
S\left( \varphi \right) =\int_{0}^{T}\left\langle \varphi \left(
X_{t}\right) ,\circ dX_{t}\right\rangle
\end{equation*}
where $\left( X_{t}\right) _{t\in \left[ 0,T\right] }$ is a stochastic
process in $\mathbb{R}^{d}$ and $\varphi $ is a $1$-form on $\mathbb{R}^{d}$%
. We give results also for It\^{o} integrals, but the Stratonovich ones do
not depend on the parametrization and therefore are more natural in a
geometrical context. Question a) amounts to ask whether the mapping $\varphi
\mapsto $ $S\left( \varphi \right) $ has a pathwise meaning and is
(pathwise) linear continuous in some topology. Question b) looks more deep,
and concernes the classification of $S$ in terms of classical types of
currents.

Concerning problem a) we deal with the following claim: 
\begin{eqnarray*}
&&\text{all semimartingales and Lyons-Zheng processes
define }1\text{-currents.}
\end{eqnarray*}
The qualitative part of this result is a by-product of known theories
even if it was not explicitly stated in the language of currents; our
contribution here is a precise estimate of the Sobolev regularity of
such currents.
The topology on $\varphi $ that we are able to consider is the $H^{s}$
topology for $s> d/2$ or $s> d/2+1$ depending on some
conditions. Our result on problem a) is essentially of analytic nature,
with little geometrical content, in contrast to problem b).

Problem a) has been considered in the literature (not with the language of
currents). Among others, two very relevant results are provided by the
integration theories of Young \cite{Young} and Lyons \cite{Lyons}. In the
theory of Young, $\left( X_{t}\right) $ is a deterministic function having
an H\"{o}lder continuous regularity with exponent $\gamma > 1/2$. For
such (and more general) curves we are even able to solve the more difficult
problem b), as described below. The outstanding (and complex) theory of
Lyons on integration of rough paths solves question a) in a variety of
cases. Restricting such theory to continuous semimartingales, it provides
pathwise integration with continuous dependence on $\varphi $ in a class
like $C^{1,\varepsilon }$ ($\varepsilon $-H\"{o}lder continuous
derivatives). Our results are not exactly comparable because of the
different topologies on $\varphi $, but we can at least say that our results
are much less deep than those of Lyons, just more straightforward.
A more abstract result which implies that semimartingales define
1-currents is a general theorem of Minlos'~\cite{Fernique} on linear
random functionals in nuclear spaces. Compared to that, the advantages
of our approach is that it gives a precise Sobolev regularity of these
currents, which is of interest in application to fluid dynamics,
see~\cite{Ascona}.  

Concerning problem b), the most common class of $1$-currents considered in
the literature on geometric measure theory are the integer multiplicity
(i.m.) rectifiable $1$-currents. This class is certainly not suitable to
incorporate the paths of the usual stochastic processes since such paths
commonly have infinite $1$-dimensional Hausdorff measure. One of the other
most important classes of currents are the \textit{integral flat chains}. We
concentrate our investigation on this class. In the 1-dimensional case, an
integral flat chain $S$ is the \textit{boundary} of an i.m. rectifiable $2$%
-current $T$, up to an additional i.m. rectifiable $1$-current $R$: 
\begin{equation*}
S\left( \varphi \right) =\partial T\left( \varphi \right) +R\left( \varphi
\right) .
\end{equation*}
In very loose terms, an integral flat $1$-chain is (a piece of) the boundary
of a quite good surface, a surface that in particular has finite $2$%
-dimensional Hausdorff measure. In intuitive geometrical terms we try to
understand whether a very irregular curve as the path of a common stochastic
process may be seen as the boundary of a surface with finite $2$-dimensional
Hausdorff measure.

The Hausdorff dimension of typical paths of Brownian motion is already $2$,
so it is very difficult to expect a positive answer in such a case (a small
room is left open by the fact that the correct gauge function for Brownian
paths is $r^{2}\log r$ in dimension $d\geq 3$ and slightly better in smaller
dimension). On the contrary, we can prove that the class of integral flat $1$%
-chains includes the paths of many processes just a little more regular than
semimartingales. We prove that 
\begin{equation*}
\text{curves of class $C^{\gamma }$ with $\gamma >1/2$ define integral flat chains.}
\end{equation*}
These classes
include for instance the typical paths of fractional Brownian motion with
Hurst parameter $H=\gamma$.

The fact that such irregular paths are boundaries of rectifiable $2$%
-currents may open interesting future investigations, like minimal surfaces
with prescribed irregular boundary, homotopy theory in such class, and
others. A direct consequence is the existence of the concept of \textit{%
integral flat distance} between two stochastic processes with paths $C^{\gamma }$, $\gamma >%
1/2$: it is a more geometrical notion
of distance  with respect to the usual analytic topologies,
independent for instance from the parametrization, and we hope it may be
useful in the future.

From the viewpoint of integration theory, our result on problem b) gives a
new definition of integral for processes with $C^{\gamma }$ paths, $\gamma >1/2$%
. By a different method such integrals have been introduced before by Young.
It is interesting to notice that Young's definition employes a telescopic
series;\ the $2$-current $T$ used here vaguely looks like a continuous
version of such series. 


For processes like Brownian motion we have not found  a classical notion
of current that may fit. The class of integral flat $1$-chains is presumably
not suitable, and at least we are sure that the surfaces we consider in
this paper have infinite $2$-dimensional Hausdorff measure for Brownian
paths, see the last section.  
However the technique we use to define the integral can be generalized 
 to provide an integration theory for
curves in classes $C^\gamma$ with $\gamma > 1/3$ which delivers
results analogous to that of T.~Lyons~\cite{Lyons}.  
These results will be reported elsewhere~\cite{GubWIP}.  


Since we address this work to different audiences, we sometimes repeat for
the reader's convenience some definition and facts that are very elementary in
one of the two theories (stochastic analysis and geometric measure theory).
We apologize with expert readers for the style that may look redundant in a
number of places.

\section{Preliminaries on deterministic and stochastic currents}

\subsection{Deterministic currents}

In this subsection we briefly recall a few definitions of geometric measure
theory used in the sequel. More informations can be found in \cite{Fed}, 
\cite{Morgan}, \cite{Simon}, \cite{GMS}.

We shall denote the Euclidean norm and scalar product in $\mathbb{R}^{d}$ ($%
d $ is fixed throughout the paper) by $\left| \cdot\right| $ and $\left\langle
\cdot,\cdot\right\rangle $ respectively. We shall also need the space $\Lambda _{2}%
\mathbb{R}^{d}$ of $2$-vectors, its dual space $\Lambda ^{2}\mathbb{R}^{d}$
of $2$-covectors, and the duality between them, still denoted by $%
\left\langle \cdot,\cdot\right\rangle $. If we represent vectors and covectors in
the standard basis (the summations are extented from $1$ to $d$) 
\begin{equation*}
v=\sum_{i<j}v^{ij}\,e_{i}\wedge e_{j},\qquad w=\sum_{i<j}w_{ij}\,e^{i}\wedge
e^{j}
\end{equation*}
then 
\begin{equation*}
\left\langle w,v\right\rangle =\sum_{i<j}w_{ij}v^{ij}.
\end{equation*}
Notationally, we prefer to write the covectors in the first argument of $%
\left\langle .,.\right\rangle $ since later on the notation for stochastic
integrals is more natural. The norm of $2$-vectors and $2$-covectors is also
defined as 
\begin{equation*}
\left| v\right| ^{2}=\sum_{i<j}\left( v^{ij}\right) ^{2},\qquad \left|
w\right| ^{2}=\sum_{i<j}\left( w_{ij}\right) ^{2}.
\end{equation*}
Notice that for all $v_{1},v_{2}\in \mathbb{R}^{d}$ and $w^{1},w^{2}\in 
\mathbb{R}^{d}$, the duality between the $2$-vector 
\begin{equation*}
v_{1}\wedge v_{2}=\sum_{i<j}\left[ \left( v_{1}\right) _{i}\left(
v_{2}\right) _{j}-\left( v_{1}\right) _{j}\left( v_{2}\right) _{i}\right]
\,e_{i}\wedge e_{j}
\end{equation*}
and the $2$-covector $w^{1}\wedge w^{2}$ similarly defined, is given by 
\begin{equation*}
\left\langle w^{1}\wedge w^{2},v_{1}\wedge v_{2}\right\rangle =\det \left(
\left\langle v_{i},w^{j}\right\rangle \right) .
\end{equation*}
Similarly we have 
\begin{equation*}
\left| v_{1}\wedge v_{2}\right| ^{2}=\det \left( \left\langle
v_{i},v_{j}\right\rangle \right) =\left| v_{1}\right| ^{2}\left|
v_{2}\right| ^{2}-\left\langle v_{1},v_{2}\right\rangle ^{2}.
\end{equation*}
Let $\mathcal{D}^{k}$ be the space of all infinitely differentiable and
compactly supported $k$-forms on $\mathbb{R}^{d}$. A $k$\textit{-dimensional
current} is a linear continuous functional on $\mathcal{D}^{k}$. We denote
by $\mathcal{D}_{k}$ the space of $k$-currents. In this paper we are only
interested in the cases $k=1,2$, so we recall a few corresponding notations.
Let $\varphi $ (resp. $\psi $) be an element of $\mathcal{D}^{1}$ (resp. $%
\mathcal{D}^{2}$). We shall write them as 
\begin{equation*}
\varphi =\sum_{i=1}^{d}\varphi _{i}\,dx^{i},\qquad \psi =\sum_{i<j}\psi
_{ij}\,dx^{i}\wedge dx^{j}
\end{equation*}
where $\varphi _{i}$ and $\psi _{ij}$ are infinitely differentiable and
compactly supported functions on $\mathbb{R}^{d}$.
Typical examples of $1$-currents are those induced by regular curves $\left(
X_{t}\right) _{t\in \left[ 0,T\right] }$ in $\mathbb{R}^{d}$: 
\begin{equation*}
\varphi \mapsto S(\varphi ):=\int_{0}^{T}\left\langle \varphi \left(
X_{t}\right) ,\dot{X}_{t}\right\rangle \,dt
\end{equation*}
while typical examples of $2$-currents are given by regular surfaces $\left(
f\left( t,s\right) \right) _{(t,s)\in A}$ in $\mathbb{R}^{d}$, where $A$ is
a Borel set in $\mathbb{R}^{2}$: 
\begin{eqnarray*}
\psi &\mapsto &T(\psi ):=\int_{A}\left\langle \psi \circ f,\frac{\partial f}{%
\partial t}\wedge \frac{\partial f}{\partial s}\right\rangle \,dtds \\
&=&\int_{A}\sum_{i<j}\left( \psi _{ij}\circ f\right) \,\left\langle
dx^{i}\wedge dx^{j},\frac{\partial f}{\partial t}\wedge \frac{\partial f}{%
\partial s}\right\rangle \,dtds \\
&=&\int_{A}\sum_{i<j}\left( \psi _{ij}\circ f\right) \,\left( \frac{\partial
f_{i}}{\partial t}\frac{\partial f_{j}}{\partial s}-\frac{\partial f_{j}}{%
\partial t}\frac{\partial f_{i}}{\partial s}\right) \,dtds.
\end{eqnarray*}
Notice that the $k$-dimensional Hausdorff measure $\mathcal{H}^{k}$ allows
one to rewrite the previous examples in a more compact and intrinsic way.
Assume for instance that $\left( X_{t}\right) _{t\in \left[ 0,T\right] }$
and $\left( f\left( t,s\right) \right) _{(t,s)\in A}$ are injective. Then 
\begin{equation*}
S(\varphi )=\int_{\mathcal{M}}\left\langle \varphi ,\xi \right\rangle \,d%
\mathcal{H}^{1},\qquad T(\psi )=\int_{\mathcal{M}}\left\langle \psi ,\xi
\right\rangle \,d\mathcal{H}^{2}
\end{equation*}
where in the first example $\mathcal{M}$ is the support of the curve $\left(
X_{t}\right) _{t\in \left[ 0,T\right] }$ and $\xi $ is the unit tangent
vector with the orientation given by $\left( X_{t}\right) _{t\in \left[ 0,T%
\right] }$, while in the second example $\mathcal{M}$ is the support of the
surface $\left( f\left( t,s\right) \right) _{(t,s)\in A}$ and $\xi $ is the
unit $2$-vector 
\begin{equation*}
\xi =\left| \frac{\partial f}{\partial t}\wedge \frac{\partial f}{\partial s}%
\right| ^{-1}\frac{\partial f}{\partial t}\wedge \frac{\partial f}{\partial s%
}
\end{equation*}
defining an orientation of the tangent manifold. Here 
\begin{equation*}
\left| \frac{\partial f}{\partial t}\wedge \frac{\partial f}{\partial s}%
\right| =\sqrt{\left| \frac{\partial f}{\partial t}\right| ^{2}\left| \frac{%
\partial f}{\partial s}\right| ^{2}-\left\langle \frac{\partial f}{\partial t%
},\frac{\partial f}{\partial s}\right\rangle ^{2}}.
\end{equation*}

These two examples can be generalized in the notion of \textit{integer
multiplicity }(\textit{i.m.} in the sequel)\textit{\ rectifiable }$k$\textit{%
-current}: it is a $k$-current $T$ of the form ($\phi \in \mathcal{D}^{k}$) 
\begin{equation*}
T\left( \phi \right) =\int_{\mathcal{M}}\left\langle \phi ,\xi \right\rangle
\theta \,d\mathcal{H}^{k}
\end{equation*}
where $\mathcal{M}$ is a $k$-rectifiable set in $\mathbb{R}^{d}$, $\theta :%
\mathcal{M}\rightarrow \mathbb{R}$ is integer valued and $\mathcal{H}%
^{k}\lfloor \mathcal{M}$ integrable, and $\xi $ is a $\mathcal{H}^{k}$%
-measurable unitary $k$-vector field on $\mathcal{M}$, with $\xi \left(
x\right) \in T_{x}\mathcal{M}$ for $\mathcal{H}^{k}\lfloor \mathcal{M}$-a.e. 
$x$ (compared to the most general definitions, here we have restricted
ourselves to sets $\mathcal{M}$ with $\mathcal{H}^{k}\left( \mathcal{M}%
\right) <\infty $). We recall that a $k$-rectifiable set $\mathcal{M}$ in $%
\mathbb{R}^{d}$ is a $\mathcal{H}^{k}$-measurable set with 
\begin{equation*}
\mathcal{H}^{k}\left( \mathcal{M}\right) <\infty
\end{equation*}
that can be decomposed as 
\begin{equation*}
\mathcal{M}=N\cup \,\bigcup_{i}f_{i}\left( A_{i}\right)
\end{equation*}
where $N$ is a set with $\mathcal{H}^{k}\left( N\right) =0$, $A_{1}$, $A_{2}$%
, ..., are Borel sets of $\mathbb{R}^{k}$, and $f_{1}$, $f_{2}$, ..., are
Lipschitz continuous maps $f_{i}:A_{i}\rightarrow \mathbb{R}^{d}$. One can
prove that there exists such a decomposition with $C^{1}$ mappings $f_{i}$
and disjoint images $f_{1}\left( A_{1}\right) $, $f_{2}\left( A_{2}\right) $%
, and so on. The $k$-rectifiable sets, although admit quite wild
singularities and intersections, have a number of properties similar to
those of $C^{1}$ $k$-manifolds. For instance, a $k$-rectifiable set $%
\mathcal{M}$ has finite $k$-dimensional Hausdorff measure; and there is a
set $N^{\prime }\subset \mathcal{M}$ with $\mathcal{H}^{k}\left( N^{\prime
}\right) =0$ such that for all $x\in \mathcal{M}\setminus N^{\prime }$ a
tangent space $T_{x}\mathcal{M}$ can be defined and a measurable orientation
can be choosen (i.e. a unitary $k$-vector field $\xi $ as above).

Finally, we recall the notion of \textit{integral flat chain}. Given
a $k$-current $T$, it is defined a $(k-1)$-current $\partial T$, the
\emph{boundary} of $T$, by the identity $\partial T\left( \phi \right) =T\left(
d\phi \right) $, where $\phi \in \mathcal{D}^{k-1}$ and $d\phi $ is its
differential. For instance, in the case $k=2$, we have 
\begin{equation*}
\phi =\sum_{i}\phi _{i}\,x^{i},\qquad d\phi =\sum_{i<j}\left( \frac{\partial
\phi _{j}}{\partial x_{i}}-\frac{\partial \phi _{i}}{\partial x_{j}}\right)
\,x^{i}\wedge x^{j}.
\end{equation*}
An integral flat $(k-1)$-chain $S$ is a $(k-1)$-current of the form 
\begin{equation*}
S\left( \phi \right) =\partial T\left( \phi \right) +R\left( \phi \right)
\end{equation*}
where $T$ is a i.m. rectifiable $k$-current and $R$ is a i.m. rectifiable $%
(k-1)$-current. In general, integral flat $(k-1)$-chain are not necessarily
i.m. rectifiable $(k-1)$-currents: the boundary of a i.m. rectifiable $k$%
-current need not to be rectifiable. These non-rectifiable boundaries are
exactly the places where we can find the irregular paths of certain
stochastic processes. Here is a trivial example of this fact for $1$%
-dimensional stochastic processes (unfortunately not very illuminating about
the difficulties arising in dimension $d>1$).

\begin{example}
Let $\left( X_{t}\right) _{t\in \left[ 0,T\right] }$ be a continuous real
valued function, for instance any typical path of a $1$-dimensional
stochastic process. Let $\lambda $ be a number smaller than the minimum of $%
\left( X_{t}\right) $. Consider the subgraph of $\left( X_{t}\right) $ above 
$\lambda $, i.e. the set 
\begin{equation*}
\mathcal{M}=\left\{ \left( t,x\right) \in \mathbb{R}^{2}:0\leq t\leq
T,\lambda \leq x\leq X_{t}\right\} .
\end{equation*}
It defines the i.m. rectifiable $2$-current in $\mathbb{R}^{2}$ given by 
\begin{equation*}
T\left( \psi \right) :=\int_{\mathcal{M}}\left\langle \psi ,e_{1}\wedge
e_{2}\right\rangle \,d\mathcal{H}^{2}=\int_{\mathcal{M}}\psi _{tx}\,d%
\mathcal{H}^{2}
\end{equation*}
where $\psi =\psi _{tx}\,dt\wedge dx$. Topologically, the boundary of $%
\mathcal{M}$ is $\left( X_{t}\right) $ plus three segments, so intuitively
we see that $\left( X_{t}\right) $ defines an integral flat $1$-chain. At
rigorous ground, denoting $1$-forms as 
\begin{equation*}
\varphi =\varphi _{t}\,dt+\varphi _{x}\,dx
\end{equation*}
the boundary of $T$ is given by 
\begin{equation*}
\partial T\left( \varphi \right) =T\left( d\varphi \right)
=\int_{0}^{T}\left( \int_{\lambda }^{X_{t}}\left( \frac{\partial \varphi
_{x}}{\partial t}-\frac{\partial \varphi _{t}}{\partial x}\right) dx\right)
dt.
\end{equation*}
Assume now that on $\left( X_{t}\right) $ we can perform rules of calculus
of classical type (as for Lipschitz functions or semimartingales by
Stratonovich stochastic calculus). Since the final result of thi example is
trivial, we do not insist on the rigorous details. Then, setting 
\begin{equation*}
\Phi _{x}\left( t,x\right) =\int_{\lambda }^{x}\varphi _{x}\left(
t,y\right) dy,
\end{equation*}
we have 
\begin{eqnarray*}
\partial T\left( \varphi \right) =\int_{0}^{T}\frac{\partial }{\partial t}%
\left[ \Phi _{x}\left( t,X_{t}\right) -\Phi _{x}\left( t,\lambda \right) %
\right] dt-\int_{0}^{T}\varphi _{x}\left( t,X_{t}\right) \circ dX_{t} \\
+\int_{0}^{T}\left( -\varphi _{t}\left( t,X_{t}\right) +\varphi _{t}\left(
t,\lambda \right) \right) dt
\end{eqnarray*}
(we have used $\circ $ to denote the suitable kind of integration required
by the previous computations) and therefore 
\begin{equation*}
\begin{split}
\int_{0}^{T}&\left\langle \varphi ,\circ \left( dt,dX_{t}\right)
\right\rangle =\int_{0}^{T}\varphi _{x}\left( t,X_{t}\right) \circ
dX_{t}+\int_{0}^{T}\varphi _{t}\left( t,X_{t}\right) dt \\
& = T\left( d\varphi \right) 
+\int_{\lambda }^{X_{T}}\varphi _{x}\left( T,y\right) dy
-\int_{\lambda
}^{X_{0}}\varphi _{x}\left( 0,y\right) dy+\int_{0}^{T}\varphi _{t}\left(
t,\lambda \right) dt.
\end{split}
\end{equation*}
Several comments are in order on this example. When $\left( X_{t}\right) $
is Brownian motion, it provides an example of i.m. rectifiable $2$-current $%
T $ with a boundary that is not an i.m. rectifiable $1$-current. It also
shows that integral flat $1$-chains may be the right objects to describe
stochastic integrals. Finally, restricted to dimension $1$ it solves
problems a) and b) posed in the introduction in complete generality.
However, with respect to problem a) this success is just a more complicate
version of the well known definition (see for instance \cite{Fo}) 
\begin{equation*}
\int_{0}^{T}\varphi _{x}\left( t,X_{t}\right) \circ dX_{t}=\Phi _{x}\left(
T,X_{T}\right) -\Phi _{x}\left( 0,X_{0}\right) -\int_{0}^{T}\frac{\partial
\Phi _{x}}{\partial t}\left( t,X_{t}\right) dt.
\end{equation*}
In other words, it is well known that in dimension one Stratonovich
stochastic integration can be defined pathwise (also It\^{o} integration can
be performed pathwise when the quadratic variation is finite, see \cite{Fo}%
). The analysis in dimension $d>1$ is completely different.
\end{example}

\subsection{Stochastic currents}

Let $\left( X_{t}\right) _{t\in \left[ 0,T\right] }$ be a continuous
semimartingale with values in $\mathbb{R}^{d}$. To each smooth $1$-form $%
\varphi $ on $\mathbb{R}^{d}$ we may associate the random variable 
\begin{equation}
S\left( \varphi \right) :=\int_{0}^{T}\left\langle \varphi \left(
X_{t}\right) ,\circ dX_{t}\right\rangle  \label{strat}
\end{equation}
where the integral is understood in the sense of Stratonovich. The mapping $%
\varphi \mapsto S\left( \varphi \right) $ is continuous with respect to the
convergence in probability. Motivated by this basic example it looks
reasonable to give the following definition.

\begin{definition}
Given a complete probability space $\left( \Omega ,\mathcal{A},P\right) $, a
stochastic $k$-current is a continuous linear mapping from the space $%
\mathcal{D}^{k}$ to the space $L^{0}\left( \Omega \right) $ of real valued
random variables on $\left( \Omega ,\mathcal{A},P\right) $, endowed with the
convergence in probability.
\end{definition}

The usual classes of stochastic processes considered in stochastic analysis
give rise to stochastic $1$-currents: semimartingales, Lyons-Zheng
processes, processes with finite $p$ variation (for suitable $p$),
fractional Brownian motion (for suitable Hurst parameter), certain Dirichlet
processes.

At this level it is difficult to see an interesting relation with classical
geometric measure theory. The link arises if we try to understand the
previous stochastic integrals in a pathwise sense. So we introduce the
following definition.

\begin{definition}
We say that the stochastic $k$-current $\varphi \mapsto S\left( \varphi
\right) $ has a \textit{pathwise realization} if there exists a measurable
mapping 
\begin{equation*}
\omega \mapsto \mathcal{S}\left( \omega \right)
\end{equation*}
from $\left( \Omega ,\mathcal{A},P\right) $ to the space $\mathcal{D}_{k}$
of deterministic currents (endowed with the natural topology of
distributions), such that 
\begin{equation}
\left[ S\left( \,\varphi \right) \right] \left( \omega \right) =\left[ 
\mathcal{S}\left( \omega \right) \right] \left( \varphi \right) \text{ \quad
for }P\text{-a.e. }\omega \in \Omega .  \label{modific}
\end{equation}
for every $\varphi \in \mathcal{D}^{k}$.
\end{definition}
In terms of these definitions we may reformulate the two problems of the
introduction as follows:
\begin{itemize}
\item[a)] given a stochastic process $\left( X_{t}\right) _{t\in \left[ 0,T\right]
} $, does there exist a pathwise realization $\mathcal{S}\left( \omega
\right) $ of the associated stochastic $1-$current $S\left( \varphi \right) $
defined by (\ref{strat})?
\item[b)] can we classify $\mathcal{S}\left( \omega \right) $ in terms of classical
currents, for $P$-a.e. $\omega $?
\end{itemize}
The existence of a pathwise realization is a difficult problem. The
difficulty is described in the remark at the beginning of the next section
in terms of selection of representatives in the equivalence classes of
stochastic integrals. In the theory of stochastic processes this is the
problem of existence of a continuous modification of a given random field $%
\left( X_{\varphi }\left( \omega \right) \right) $. Here the parameter of
the field is $\varphi \in \mathcal{D}^{k}$, so the parameter space is
infinite dimensional and well-known criteria like the Kolmogorov regularity
theorem do not apply (some generalizations are known in the literature but
their effective use is very limited). The problem of existence of a
continuous modification of a random field with infinite dimensional
parameter space has been studied and some general ideas have been developed,
but usually it is better to find out ad hoc methods, as we shall do. We just
recall now two general criteria from \cite{Sk} (see also \cite{Fbook}).
\begin{lemma}
Let $\varphi \mapsto S\left( \varphi \right) $ be a linear continuous
mapping from a separable Banach space $E$ to $L^{0}\left( \Omega \right) $
(with the convergence in probability). Assume that there exists a random
variable $C\left( \omega \right) $ such that for all given $\varphi \in E$
we have 
\begin{equation*}
\left| S\left( \varphi \right) \left( \omega \right) \right| \leq C\left(
\omega \right) \left\| \varphi \right\| _{E}\quad \text{for }P\text{-a.e. }%
\omega \in \Omega .
\end{equation*}
Then there exists a measurable mapping $\omega \mapsto $ $\mathcal{S}\left(
\omega \right) $ from $\left( \Omega ,\mathcal{A},P\right) $ to the dual $%
E^{\prime }$ such that for all given $\varphi \in E$ we have (\ref{modific})
(hence $\mathcal{S}\left( \omega \right) $ is a pathwise realization of $%
S\left( \varphi \right) $).
\end{lemma}
We shall use this criterium in the next section. It is not very powerful
since its assumption is a pathwise estimate (it is almost a tautology). The
proof is elementary and can be found in the above mentioned references. More
interesting is the following criterium since it is based on an assumption in
mean square. However, in the next section we shall not use it directly but
an ad hoc argument based on Fourier transform, that looks more flexible. We
skech the proof (contained in \cite{Sk}, \cite{Fbook}), for comparison with
the method of the next section.
\begin{lemma}
Let $\varphi \mapsto S\left( \varphi \right) $ be a linear continuous
mapping from a separable Hilbert space $H$ to $L^{2}\left( \Omega \right) $.
Assume that it is Hilbert-Schmidt: for some complete orthonormal system $%
\left\{ e_{i}\right\} $ in $H$ we have 
\begin{equation*}
\sum_{i=1}^{\infty }\left| S\left( e_{i}\right) \right| _{L^{2}\left( \Omega
\right) }^{2}<\infty .
\end{equation*}
Then there exists a measurable mapping $\omega \mapsto $ $\mathcal{S}\left(
\omega \right) $ from $\left( \Omega ,\mathcal{A},P\right) $ to $H$, i.e. a
random vector of $H$, such that for all given $\varphi \in H$ we have (\ref
{modific})
\end{lemma}
\begin{proof}
Schwartz inequality gives us 
\begin{eqnarray*}
\left| S\left( \varphi \right) \left( \omega \right) \right| ^{2}=\left|
\sum_{i=1}^{\infty }S\left( e_{i}\right) \left( \omega \right) \left\langle
\varphi ,e_{i}\right\rangle _{H}\right| ^{2} \\
\leq \sum_{i=1}^{\infty }\left| S\left( e_{i}\right) \left( \omega \right)
\right| ^{2}\sum_{i=1}^{\infty }\left\langle \varphi ,e_{i}\right\rangle
_{H}^{2} \\
=C\left( \omega \right) \left\| \varphi \right\| _{H}^{2}.
\end{eqnarray*}
The non negative (a priori possibly infinite) r.v. $C\left( \omega \right) $
has finite mean by assumption, then it finite a.s. So we may apply the first
lemma.
\end{proof}

\section{Paths of semimartingales and Lyons-Zheng processes define $1$%
-currents}

\subsection{The case of semimartingales}

In the sequel we tacitly assume that processes are defined on a complete
probability space $\left( \Omega ,\mathcal{A},P\right) $, with expectation
denoted by $E$. We also assume that a standard filtration $\mathcal{F}%
=\left( \mathcal{F}_{t}\right) $ is given, so\ that concepts like martingale
or adaptedness are referred to this filtration.

For the definition of semimartingale and corresponding integrals, see \cite
{RY}, \cite{Kunita} or many other references. We just recall a few facts
directly used below. A continuous semimartingale $\left( X_{t}\right) _{t\in %
\left[ 0,T\right] }$ is the sum of a continuous local martingale $\left(
M_{t}\right) $ and a continuous adapted process of bounded variation $\left(
V_{t}\right) $. The decomposition is unique. Given a continuous adapted
process $\left( Y_{t}\right) $ in $\mathbb{R}^{d}$, the It\^{o} integral $%
\int_{0}^{T}\left\langle Y_{t},dX_{t}\right\rangle $ is defined as 
\begin{equation*}
\int_{0}^{T}\left\langle Y_{t},dX_{t}\right\rangle :=\int_{0}^{T}\left\langle
Y_{t},dM_{t}\right\rangle +\int_{0}^{T}\left\langle Y_{t},dV_{t}\right\rangle
\end{equation*}
and similarly\ for the Stratonovich integral $\int_{0}^{T}\left\langle
Y_{t},\circ dX_{t}\right\rangle $, where now $\left( Y_{t}\right) $ is
assumed to be either a continuous semimartingale or $Y_{t}=\varphi \left(
X_{t}\right) $ with the $1$-form $\varphi $ is of class $C^{1}$ (one can
unify these two cases with the language of Dirichlet processes). The
previous integrations in $dV_{t}$ are classical pathwise integrations in the
Riemann-Stieltjes sense, while the integrals with respect to $\left(
M_{t}\right) $ are the following limits in probability (they exist under the
previous assumptions on $\left( Y_{t}\right) $): 
\begin{equation*}
\int_{0}^{T}\left\langle Y_{t},dM_{t}\right\rangle := \operatornamewithlimits{P-lim}_{n\rightarrow
\infty }\sum_{t_{i}\in \pi _{n}}\left\langle
Y_{t_{i}},M_{t_{i+1}}-M_{t_{i}}\right\rangle
\end{equation*}
\begin{equation*}
\int_{0}^{T}\left\langle Y_{t},\circ dM_{t}\right\rangle :=
\operatornamewithlimits{P-lim}_{n\rightarrow \infty }\sum_{t_{i}\in \pi _{n}}\left\langle \frac{%
Y_{t_{i+1}}+Y_{t_{i}}}{2},M_{t_{i+1}}-M_{t_{i}}\right\rangle
\end{equation*}
where $\pi _{n}$ is any sequence of partitions of $\left[ 0,T\right] $
converging to zero.

\begin{remark}
These integrals are $P$-equivalence classes. The evaluation at a given $%
\omega $, namely the pathwise integration, is a priori meaningless. Given $%
\varphi $ one may of course take a representative in the equivalence class
and have a meaning for all $\omega $, but an arbitrarity choice of the
representative cannot give us any good property (even the linearity) of the
mapping $\varphi \mapsto S\left( \varphi \right) $ evaluated at single
points $\omega $. The existence of a pathwise realization means that it is
possible to choose representatives in such a way that the mapping $\varphi
\mapsto S\left( \varphi \right) $, evaluated at almost every given point $%
\omega $, is linear and  continuous.
\end{remark}

\bigskip Given two continuous semimartingales $\left( X_{t}\right) $ and $%
\left( Y_{t}\right) $ in $\mathbb{R}^{d}$, one can define the quadratic 
\textit{covariation} processes between their components $\left[ X^{\alpha
},Y^{\beta }\right] _{t}$ as 
\begin{equation*}
\left[ X^{\alpha },Y^{\beta }\right] _{t} := \operatornamewithlimits{P-lim}_{n\rightarrow \infty }\sum 
_{\substack{ t_{i}\in \pi _{n}  \\ t_{i}\leq t}}\left\langle
X_{t_{i+1}}^{\alpha }-X_{t_{i}}^{\alpha },Y_{t_{i+1}}^{\beta
}-Y_{t_{i}}^{\beta }\right\rangle .
\end{equation*}
One has $[ X^{\alpha },Y^{\beta }] _{t}=[ M_{X}^{\alpha
},M_{Y}^{\beta }] _{t}$ where $M_{X}$ and $M_{Y}$ are the martingale
parts of $\left( X_{t}\right) $ and $\left( Y_{t}\right) $ (the bounded
variation terms do not contribute to the quadratic variation). The relation
between Stratonovich and It\^{o} integral is now 
\begin{equation*}
\int_{0}^{T}\left\langle Y_{t},\circ dX_{t}\right\rangle
=\int_{0}^{T}\left\langle Y_{t},dX_{t}\right\rangle +1/2\sum_{\alpha
=1}^{d}\left[ Y^{\alpha },X^{\alpha }\right] _{T}
\end{equation*}
as one may easily check by means of the finite sums. Similar facts hold true
when $Y_{t}=\varphi \left( X_{t}\right) $, with the additional formula: 
\begin{equation*}
\left[ X^{\alpha },\varphi _{\beta }\left( X\right) \right]
_{T}=\int_{0}^{T}\sum_{\delta=1}^{d}\frac{\partial \varphi _{\beta }}{\partial
x_{\delta}}\left( X_{t}\right) d\left[ M^{\delta},M^{\alpha }\right] _{t}
\end{equation*}
where $\left( M_{t}\right) $ is the martingale part of $\left( X_{t}\right) $%
.

Finally, we recall the Burkholder-Davis-Gundy inequality. For all $p\geq 1$,
there exists a constant $C_{p}>0$ such that 
\begin{equation*}
E\left[ \left| \int_{0}^{T}Y_{t}dM_{t}^{\delta}\right| ^{2p}\right] \leq C_{p}E%
\left[ \left| \int_{0}^{T}\left| Y_{t}\right| ^{2}d\left[ M^{\delta}\right]
_{t}\right| ^{p}\right]
\end{equation*}
where $\left( Y_{t}\right) $ is any continuous adapted scalar process.

Let us come to our results. Denote the Fourier transform of $\varphi \left(
x\right) $ by $\hat{\varphi}\left( k\right) $: 
\begin{equation*}
\hat{\varphi}\left( k\right) :=\int_{\mathbb{R}^{d}}e^{-i\left\langle
k,x\right\rangle }\varphi \left( x\right) \,dx.
\end{equation*}
The following simple lemma will be our key ingredient and we guess it may be
useful in other contexts. It is inspired by the vision of stochastic
integrals as currents (i.e. as generalized random fields, so that it is not
strange to perform their Fourier transform) and by the computations of \cite
{FGub}.
\begin{lemma}
Let $\left( M_{t}\right) $ be an $L^{2}$-bounded continuous martingale ($%
E[ | M_{t}| ^{2}]
< \infty$ for $t \in [0,T]$), $(X_{t})$ be a continuous adapted process, and $\varphi $ be in $%
\mathcal{D}^{1}$. Then 
\begin{equation}
\int_{0}^{T}\left\langle \varphi \left( X_{t}\right) ,dM_{t}\right\rangle
=\int_{\mathbb{R}^{d}}\left\langle \hat{\varphi}\left( k\right)
,Z_{k}\right\rangle \,dk\quad P\text{-a.s.}  \label{fouriercharact}
\end{equation}
where 
\begin{equation*}
Z_{k}:=\int_{0}^{T}e^{i\left\langle k,X_{t}\right\rangle }dM_{t}.
\end{equation*}
A similar result holds true for $\left( V_{t}\right) $ in place of $\left(
M_{t}\right) $, when $\left\| V\right\| _{var}^{2}\in L^{1}\left( \Omega
\right) $. Moreover, a similar result holds true for the Stratonovich
integral when $(X_{t})$ is a semimartingale where the martingale and bounded
variation parts satisfy the same integrability assumptions of $\left(
M_{t}\right) $ and $\left( V_{t}\right) $ (so in particular for $X=M+V$).
\end{lemma}
\begin{proof}
We give the proof only in the first case, since the others are entirely
similar. First notice that the mapping $(\omega ,k)\mapsto Z_{k}\left(
\omega \right) $ is measurable (by the formula for the It\^{o} integral as
limit of finite sums) and the function $\left\langle \hat{\varphi}\left(
k\right) ,Z_{k}\right\rangle $ is jointly integrable in $(\omega ,k)$, so
the right-hand-side of (\ref{fouriercharact}) is well defined. Indeed 
\begin{equation}
E[ \left| Z_{k}\right| ^{2}] \leq C_{d}E\int_{0}^{T}d\left[ M%
\right] _{t}=C_{d}E\left[ M\right] _{T}<\infty  \label{easybound}
\end{equation}
(the last inequality is due to the assumption on $\left( M_{t}\right) $ and
Corollary 1.25, Ch. IV of \cite{RY}) so $E \left| Z_{k}\right| 
\leq (E\left| Z_{k}\right| ^{2})^{1/2}\leq (C_{d}E\left[ M%
\right] _{T})^{1/2}$ and therefore 
\begin{eqnarray*}
E\left[ \int_{\mathbb{R}^{d}}\left| \left\langle \hat{\varphi}\left(
k\right) ,Z_{k}\right\rangle \right| \,dk\right] \leq \int_{\mathbb{R}%
^{d}}\left| \hat{\varphi}\left( k\right) \right| E [\left| Z_{k}\right| ]%
 \,dk \\
\leq \sqrt{C_{d}E\left[ M\right] _{T}\text{ }}\int_{\mathbb{R}^{d}}\left| 
\hat{\varphi}\left( k\right) \right| \,dk
\end{eqnarray*}
where the last integral converges because of the decay properties of $\left| 
\hat{\varphi}\left( k\right) \right| $.

Formally the result (\ref{fouriercharact}) is a consequence of the heuristic
formula 
\begin{equation*}
\int_{0}^{T}\left\langle \varphi \left( X_{t}\right) ,dM_{t}\right\rangle
=\int_{\mathbb{R}^{d}}\left\langle \varphi \left( x\right) ,f\left( x\right)
\right\rangle \,dx
\end{equation*}
where 
\begin{equation*}
f\left( x\right) :=\int_{0}^{T}\delta \left( x-X_{t}\right) \,dM_{t}.
\end{equation*}
Let $p_{\varepsilon }\left( x\right) $ denote the heat kernel $\left(
2\pi \varepsilon \right) ^{-d/2}\exp ( -|x| ^{2}/(2\varepsilon)) $ and let 
\begin{equation*}
f_{\varepsilon }\left( x\right) :=\int_{0}^{T}p_{\varepsilon }\left(
x-X_{t}\right) \,dM_{t}.
\end{equation*}
Then, by stochastic Fubini theorem, \cite{RY} p. 167, 
\begin{equation*}
\int_{\mathbb{R}^{d}}\left\langle \varphi \left( x\right) ,f_{\varepsilon
}\left( x\right) \right\rangle \,dx=\int_{0}^{T}\left\langle \int_{\mathbb{R}%
^{d}}\varphi \left( x\right) p_{\varepsilon }\left( x-X_{t}\right)
\,\,dx,\,dM_{t}\right\rangle .
\end{equation*}
Hence, by Parseval theorem (we exchange the order for comparison with (\ref
{fouriercharact}) to be proved) 
\begin{equation}
\int_{0}^{T}\left\langle \varphi _{\varepsilon }\left( X_{t}\right)
,\,dM_{t}\right\rangle =\int_{\mathbb{R}^{d}}\left\langle \hat{\varphi}%
\left( k\right) ,\hat{f}_{\varepsilon }\left( k\right) \right\rangle \,dk
\label{fouriercharactappr}
\end{equation}
where $\varphi _{\varepsilon }=p_{\varepsilon }\ast \varphi $. We also have 
\begin{equation*}
\hat{f}_{\varepsilon }\left( k\right) =\int_{0}^{T}e^{i\left\langle
k,X_{t}\right\rangle }\hat{p}_{\varepsilon }\left( k\right) \,dM_{t}.
\end{equation*}
Since $\hat{p}_{\varepsilon }\left( k\right) \rightarrow 1$ uniformly on
compact sets of $k$, as $\varepsilon \rightarrow 0$, and $\varphi
_{\varepsilon }\left( X_{t}\right) $ converges to $\varphi \left(
X_{t}\right) $ uniformly in $t$, $P$-a.s., and the convergence is dominated
by a constant, we have that $\hat{f}_{\varepsilon }\left( k\right) $ and $%
\int_{0}^{T}\left\langle \varphi _{\varepsilon }\left( X_{t}\right)
,\,dM_{t}\right\rangle $ converge to\ $Z_{k}$ for all $k$ and to $%
\int_{0}^{T}\left\langle \varphi \left( X_{t}\right) ,\,dM_{t}\right\rangle $
respectively, in mean square. Therefore, first, the l.h.s. of (\ref
{fouriercharactappr}) converges to the one of (\ref{fouriercharact})\ in
mean square. As to the r.h.s., 
\begin{eqnarray*}
E\left[ \int_{\mathbb{R}^{d}}\left| \left\langle \hat{\varphi}\left(
k\right) ,\hat{f}_{\varepsilon }\left( k\right) -Z_{k}\right\rangle \right|
\,dk\right] \leq \int_{\mathbb{R}^{d}}\left| \hat{\varphi}\left( k\right)
\right| E\left[ \left| \hat{f}_{\varepsilon }\left( k\right) -Z_{k}\right| %
\right] \,dk \\
\leq \int_{\mathbb{R}^{d}}\left| \hat{\varphi}\left( k\right) \right| \,%
\sqrt{C_{d}E\left[ \left| \hat{f}_{\varepsilon }\left( k\right)
-Z_{k}\right| ^{2}\right] \text{ }}dk.
\end{eqnarray*}
The term $E[| \hat{f}_{\varepsilon }( k) -Z_{k}|^{2}] $ converges to zero for every $k$, and is bounded by a constant
(it is easily proved as (\ref{easybound})). From the decay properties of $%
\left| \hat{\varphi}\left( k\right) \right| $ we deduce that $\langle 
\hat{\varphi}( k) ,\hat{f}_{\varepsilon }( k)
-Z_{k}\rangle $ converges to zero in $L^{1}$ with respect to $(\omega
,k)$ so $\int_{\mathbb{R}^{d}}\langle \hat{\varphi}( k) ,%
\hat{f}_{\varepsilon }( k) \rangle \,dk$ converges to $%
\int_{\mathbb{R}^{d}}\langle \hat{\varphi}( k)
,Z_{k}\rangle \,dk$ in $L^{1}$ with respect to $\omega $. This
completes the proof of (\ref{fouriercharact}).
\end{proof}

\begin{theorem}
Let $\left( X_{t}\right) $ be a semimartingale in $\mathbb{R}^{d}$ of the
form $X_{t}=M_{t}+V_{t}$ as above. Consider the stochastic $1-$current $%
S\left( \varphi \right) $ defined by the Stratonovich integral 
\begin{equation*}
S\left( \varphi \right) =\int_{0}^{T}\left\langle \varphi \left(
X_{t}\right) ,\circ dX_{t}\right\rangle
\end{equation*}
and the stochastic $1-$current $I\left( \varphi \right) $ defined by the
It\^{o} integral 
\begin{equation*}
I\left( \varphi \right) =\int_{0}^{T}\left\langle \varphi \left(
X_{t}\right) ,dX_{t}\right\rangle .
\end{equation*}
Then $\varphi \mapsto S\left( \varphi \right) $ has a pathwise realization $%
\mathcal{S} $, with 
\begin{equation*}
\mathcal{S}\left( \omega \right) \in H^{-s-1}\left( \mathbb{R}^{d},\mathbb{R}%
^{d}\right) \text{ \qquad }P\text{-a.s.}
\end{equation*}
for all $s>\frac{d}{2}$, and $\varphi \mapsto I\left( \varphi \right) $ has
a pathwise realization $\mathcal{I} $, with 
\begin{equation*}
\mathcal{I}\left( \omega \right) \in H^{-s}\left( \mathbb{R}^{d},\mathbb{R}%
^{d}\right) \text{ \qquad }P\text{-a.s.}
\end{equation*}
If in addition 
\begin{equation}
\left[ M^{j},M^{i}\right] \equiv 0\text{ for }i\neq j\text{ and }\left[ M^{i}%
\right] =m_{t}\text{ for all }i  \label{algass}
\end{equation}
and for some increasing process $\left( m_{t}\right) $, then 
\begin{equation*}
\mathcal{S}\left( \omega \right) \in H^{-s}\left( \mathbb{R}^{d},\mathbb{R}%
^{d}\right) \text{ \qquad }P\text{-a.s.}
\end{equation*}
(the same result holds true for reversible semimartingales, see the next
theorem). Moreover, if $\left( M_{t}\right) $ is a square integrable
martingale and $\left\| V\right\| _{var}^{2}\in L^{1}\left( \Omega \right) $%
, then 
\begin{equation*}
\mathcal{S}\left( .\right) \in L^{2}\left( \Omega ,H^{-s-1}\left( \mathbb{R}%
^{d},\mathbb{R}^{d}\right) \right)
\end{equation*}
\begin{equation*}
\mathcal{I}\left( .\right) \in L^{2}\left( \Omega ,H^{-s}\left( \mathbb{R}%
^{d},\mathbb{R}^{d}\right) \right)
\end{equation*}
and under the assumption (\ref{algass}) 
\begin{equation*}
\mathcal{S}\left( .\right) \in L^{2}\left( \Omega ,H^{-s}\left( \mathbb{R}%
^{d},\mathbb{R}^{d}\right) \right) .
\end{equation*}
Finally, except for the result under assumption (\ref{algass}), the same
results hold true for the It\^{o} integral 
\begin{equation*}
\tilde{I}\left( \varphi \right) =\int_{0}^{T}\langle \varphi ( 
\tilde{X}_{t}) ,dX_{t}\rangle
\end{equation*}
and the analogous Stratonovich integral, when $( \tilde{X}_{t}) $
is another semimartingale in $\mathbb{R}^{d}$ (with integrability
assumptions similar to those of $( X_{t}) $ for the last results
on summability).
\end{theorem}
\begin{proof}
\textbf{Step 1} (localized problem and basic estimates; It\^{o} integral).
Let $\tau _{n}^{M}$ be a sequence of stopping times that localizes $\left(
M_{t}\right) $. Let $\tau _{n}^{\left[ M\right] }$ be the one defined as 
\begin{equation*}
\tau _{n}^{\left[ M\right] }=\inf \left\{ t\geq 0:\left[ M\right] _{t}\geq
n\right\}
\end{equation*}
when this set is non empty, $\tau _{n}^{\left[ M\right] }=T$ otherwise.
Similarly, let $\tau _{n}^{V}$ be defined as 
\begin{equation*}
\tau _{n}^{V}=\inf \left\{ t\geq 0:\left\| V\right\| \geq n\right\} .
\end{equation*}
when this set is non empty, $\tau _{n}^{V}=T$ otherwise. Finally, let $\tau
_{n}$ be defined as $\tau _{n}=\tau _{n}^{M}\wedge \tau _{n}^{\left[ M\right]
}\wedge \tau _{n}^{V}$. It localizes $\left( M_{t}\right) $ by Doob's stopping
theorem, so $\left( M_{t}^{\left( n\right) }\right) $ defined as 
\begin{equation*}
M_{t}^{\left( n\right) }=M_{t\wedge \tau _{n}}
\end{equation*}
is a martingale, and in addition $\left[ M^{\left( n\right) }\right]
_{t}\leq n$. Moreover, setting $V_{t}^{\left( n\right) }:=V_{t\wedge \tau
_{n}}$, we have $\left\| V\right\| \leq n$. Let us set $X_{t}^{\left(
n\right) }=M_{t}^{\left( n\right) }+V_{t}^{\left( n\right) }$ and introduce
the stochastic current 
\begin{equation*}
I_{n}\left( \varphi \right) :=\int_{0}^{T}\left\langle \varphi (
X_{t}^{\left( n\right) }) ,dX_{t}^{\left( n\right) }\right\rangle .
\end{equation*}
By the previous lemma we have 
\begin{equation*}
I_{n}\left( \varphi \right) =\int_{\mathbb{R}^{d}} \langle \hat{\varphi}%
( k ) ,Z_{k}^{( n) }\rangle \,dk
\end{equation*}
where 
\begin{equation*}
Z_{k}^{\left( n\right) }:=\int_{0}^{T}e^{i\langle k,X_{t}^{(
n) }\rangle }dX_{t}^{( n) }.
\end{equation*}
On these ``Fourier coefficients'' we have the estimate 
\begin{equation*}
\begin{split}
E\left[ \left| Z_{k}^{\left( n\right) }\,\right| ^{2}\right] & \leq 2E\left[
\left| \int_{0}^{T}e^{i\langle k,X_{t}^{( n) }\rangle
}dM_{t}^{\left( n\right) }\,\right| ^{2}\right] +2E\left[ \left|
\int_{0}^{T}e^{i\langle k,X_{t}^{( n) }\rangle
}dV_{t}^{\left( n\right) }\,\right| ^{2}\right] \\
& \leq C_{d}E\int_{0}^{T}d\left[ M^{\left( n\right) }\right] _{t}+C_{d}E\left[
\left( \int_{0}^{T}d| V^{(n) }| \,_{t}\right) ^{2}%
\right] \\
& =C_{d}E\left[ M^{\left( n\right) }\right] _{T}+C_{d}E\left| V^{\left(
n\right) }\right| _{T}\leq 2C_{d}n.
\end{split}
\end{equation*}
We now have 
\begin{equation*}
\begin{split}
\left| I_{n}\left( \varphi \right) \right| 
& \leq \left( \int_{\mathbb{R}^{d}}\frac{| Z_{k}^{( n) }|^{2}}{(1+|k|^{2})^{s}}\,\,dk\right) ^{\frac{%
1}{2}}\left( \int_{\mathbb{R}^{d}}\,| \hat{\varphi}( k)|^{2}
( 1+| k|^{2})^{s}\,dk\right)^{1/2} \\
& \leq C\left( \omega \right) \left\| \varphi \right\| _{H^{s}}
\end{split}
\end{equation*}
where 
\begin{equation*}
\begin{split}
E\left[ \left| C\right| ^{2}\right] & =E\int_{\mathbb{R}^{d}}
\frac{|Z_{k}^{( n) }|^{2}}{( 1+| k|^{2}) ^{s}}\,\,dk 
\leq 2C_{d}n\int_{\mathbb{R}^{d}}\frac{1}{( 1+| k|
^{2}) ^{s}}\,dk<\infty
\end{split}
\end{equation*}
for $s>d/2$. Therefore $C\left( \omega \right) $ is finite $P$-a.s.
and the lemma above applies with $E=H^{s}\left( \mathbb{R}^{d},\mathbb{R}%
^{d}\right) $. We have proved that the stochastic current $\varphi \mapsto
I_{n}\left( \varphi \right) $ has a pathwise realization $\mathcal{I}%
_{n}\left( \omega \right) $.

\textbf{Step 2} (Stratonovich integral, general case). Consider now the
stochastic current $\varphi \mapsto S_{n}\left( \varphi \right) $ defined as 
\begin{equation*}
S_{n}\left( \varphi \right) :=\int_{0}^{T}\langle \varphi (
X_{t}^{( n) }) ,\circ dX_{t}^{( n)
}\rangle .
\end{equation*}
By the relation between Stratonovich and It\^{o} integrals it is intuitively
clear that we should have the same result with one more derivative of $%
\varphi $, i.e. the topology $H^{s+1}\left( \mathbb{R}^{d},\mathbb{R}%
^{d}\right) $ on $\varphi $. Let us prove the result. We have 
\begin{equation*}
S_{n}\left( \varphi \right) =\int_{\mathbb{R}^{d}}\langle \hat{\varphi}%
\left( k\right) ,Z_{k}^{\left( n\right) }\rangle \,dk
\end{equation*}
where now we set 
\begin{equation*}
Z_{k}^{\left( n\right) }:=\int_{0}^{T}e^{i\langle k,X_{t}^{(
n) }\rangle }\circ dX_{t}^{( n) }.
\end{equation*}
We have 
\begin{equation*}
Z_{k}^{\left( n\right) \beta }=\int_{0}^{T}e^{i\langle k,X_{t}^{(
n) }\rangle }dX_{t}^{( n) \beta }+\frac{i}{2}%
\sum_{\alpha =1}^{d}k_{\alpha }\int_{0}^{T}e^{i\langle k,X_{t}^{(
n) }\rangle }d[ X^{( n) \alpha },X^{(
n) \beta }] _{t}
\end{equation*}
so that, by the estimate of the previous step, 
\begin{equation*}
E\left[ \left| Z_{k}^{\left( n\right) }\,\right| ^{2}\right] \leq
4C_{d}n+C_{d}^{\prime }\left| k\right| ^{2}\sum_{\alpha ,\beta =1}^{d}E\left[
\left| \int_{0}^{T}e^{i\langle k,X_{t}^{( n) }\rangle
}d[ X^{( n) \alpha },X^{( n) \beta }]
_{t}\right| ^{2}\right] .
\end{equation*}
Recall that $[ X^{\left( n\right) \alpha },X^{\left( n\right) \beta }%
] _{t}$ is a bounded variation function and is given by 
\begin{equation*}
\frac{1}{4}[ X^{\left( n\right) \alpha }+X^{\left( n\right) \beta }%
] _{t}-\frac{1}{4}[ X^{\left( n\right) \alpha }-X^{\left(
n\right) \beta }] _{t}
\end{equation*}
that provides the decomposition as difference of non decreasing functions.
Each of them can be controlled by $[ X^{\left( n\right) \alpha }]
_{t}$ and $[ X^{\left( n\right) \beta }] _{t}$ . Therefore we
finally have, for a new constant, 
\begin{equation}
\label{eq:strat_est1}
E\left[ | Z_{k}^{\left( n\right) }| ^{2}\right] \leq
C_{d}n ( 1+\left| k\right| ^{2}) .
\end{equation}
Repeating the argument of the previous step we get 
\begin{equation*}
\left| S_{n}\left( \varphi \right) \right| \leq C\left( \omega \right)
\left\| \varphi \right\| _{H^{s+1}}
\end{equation*}
with $C\left( \omega \right) $ a.s. finite. This proves that $\varphi
\mapsto S_{n}\left( \varphi \right) $ has a pathwise realization $\mathcal{S}%
_{n}\left( \omega \right) $, continuous in the $H^{s+1}$-topology.

\textbf{Step 3} (Stratonovich integral, under assumption (\ref{algass}%
)). With the notations of step 2, where we drop $n$ for simplicity of
notations, we decompose $Z_{k}$ in the direction of $k$ and its orthogonal
by means of a suitable projection $p_{k}$: 
\begin{equation*}
Z_{k}=\frac{k}{\left| k\right| ^{2}}\left\langle Z_{k},k\right\rangle
+p_{k}Z_{k}.
\end{equation*}
We have (using It\^{o} formula in the first line \cite{RY},\cite{Kunita},
and the relation between Stratonovich and It\^{o} integrals in the second
one) 
\begin{equation*}
\left\langle Z_{k},k\right\rangle =\int_{0}^{T }e^{i\left\langle
k,X_{t}\right\rangle }\circ d\left\langle k,X_{t}\right\rangle
=-ie^{i\left\langle k,X_{T}\right\rangle }+ie^{i\left\langle
k,X_{0}\right\rangle }
\end{equation*}
\begin{equation*}
Z_{k}^{\beta }=\int_{0}^{T }e^{ik\cdot X_{t}}dX_{t}^{\beta }+\frac{i%
k_{\beta }}{2}\int_{0}^{T }e^{ik\cdot X_{t}}\,dm_{t}
\end{equation*}
by the assumption on the covariation, and therefore 
\begin{equation*}
p_{k}Z_{k}=\int_{0}^{T }e^{ik\cdot X_{t}}d\left( p_{k}X_{t}\right) .
\end{equation*}
Summarizing we have 
\begin{equation*}
  \begin{split}
Z_{k} & 
= \frac{-ik}{\left| k\right| ^{2}}\left( e^{i\left\langle
k,X_{T}\right\rangle }-e^{i\left\langle k,X_{0}\right\rangle }\right)
+\int_{0}^{T }e^{ik\cdot X_{t}}d\left( p_{k}M_{t}\right) +\int_{0}^{T
}e^{ik\cdot X_{t}}d\left( p_{k}V_{t}\right) .    
  \end{split}
\end{equation*}
It is now easy, with estimates similar to those above, to prove that 
\begin{equation*}
E\left[ \left| Z_{k}\,\right| ^{2}\right] \leq C_{d}n.
\end{equation*}
Notice that under the additional assumption (\ref{algass}) we do not have
the factor $( 1+\left| k\right| ^{2}) $ which appears in the
estimate~(\ref{eq:strat_est1}).  Therefore $\varphi
\mapsto S_{n}\left( \varphi \right) $ has a pathwise realization $\mathcal{S}%
_{n}\left( \omega \right) $, continuous in the $H^{s}$-topology.

\textbf{Step 4} (conclusion).  Consider for instance the case of the
Stratonovich integral (the other is similar). By the locality property of
stochastic integrals, \cite{RY} proposition 2.11 of Ch. IV, on the event $%
\Omega _{n}=\left\{ \tau _{n}=T\right\} $ we have $S_{n}\left( \varphi
\right) =S\left( \varphi \right) $ a.s., for any given $\varphi $. The
sequence $\left( \Omega _{n}\right) $ increases to $\Omega $. Define $%
\mathcal{S}\left( \omega \right) $ as 
\begin{equation*}
\mathcal{S}\left( \omega \right) =\mathcal{S}_{n}\left( \omega \right) \text{
for all }\omega \in \Omega _{n}\text{.}
\end{equation*}
The definition is a.s. correct because $\mathcal{S}_{n+1}\left( \omega
\right) =\mathcal{S}_{n}\left( \omega \right) $ a.s. on $\Omega _{n}$.
Indeed, given $\varphi $, we have 
\begin{equation*}
\mathcal{S}_{n+1}\left( \omega \right) \varphi =S_{n+1}\left( \varphi
\right) \left( \omega \right) \text{ for a.e. }\omega \in \Omega
\end{equation*}
hence by the locality property 
\begin{equation*}
\mathcal{S}_{n+1}\left( \omega \right) \varphi =S\left( \varphi \right)
\left( \omega \right) \text{ for a.e. }\omega \in \Omega _{n+1}
\end{equation*}
and similarly 
\begin{equation*}
\mathcal{S}_{n}\left( \omega \right) \varphi =S\left( \varphi \right) \left(
\omega \right) \text{ for a.e. }\omega \in \Omega _{n}.
\end{equation*}
This proves that $\mathcal{S}_{n+1}\left( \omega \right) =\mathcal{S}%
_{n}\left( \omega \right) $ a.s. on $\Omega _{n}$, that the definition is
correct and that 
\begin{equation*}
\mathcal{S}\left( \omega \right) \varphi =S\left( \varphi \right) \left(
\omega \right) \text{ for a.e. }\omega \in \Omega .
\end{equation*}
Therefore we have found a pathwise realization of $S$. Finally, the
estimates in mean value and the generalization to $( \tilde{X}%
_{t}) $ can be obtained just by inspection in the previous arguments
and inequalities.
\end{proof}

\begin{remark}
The previous result and proof is very related to \cite{FGub}, although the
aim is different. The proof is also related to the one of the last lemma of
the previous section.
\end{remark}

\begin{remark}
Condition (\ref{algass}) is fulfilled for instance by the $d$-dimensional
Brownian motion $\left( W_{t}\right) $ since $\left[ W^{i},W^{j}\right]
_{t}=\delta _{ij}t$.
\end{remark}

\begin{remark}
The result under assumption (\ref{algass}) is slightly surprising. Indeed,
recall that $H^{s}$ is embedded into the space of continuous $1$-forms, and $%
H^{s+1}$ into the continuously differentiable ones. Therefore the results in
the general cases have a correspondence with the fact that It\^{o} integrals
are defined when $\varphi $ is continuous and Stratonovich one when $\varphi 
$ is continuously differentiable. But the result under assumption (\ref
{algass}) says that we have a well defined Stratonovich integral, even
pathwise defined, for all functions $\varphi \in H^{s}$, that are not
necessarily continuously differentiable. This result should be compared more
carefully with results on stochastic integration obtained when the processes
have densities with respect to the Lebesgue measure. We do not stress this
direction here.
\end{remark}

\subsection{The case of process with Lyons-Zheng structure}

The concept of process with Lyons-Zheng structure has been introduced and
studied by \cite{LZ}, \cite{LZ2}, \cite{RV}, among other references. We
follow the presentation of \cite{RV}. We say that $\left( X_{t}\right)
_{t\in \left[ 0,T\right] }$ is a Lyons-Zheng process if it has the form 
\begin{equation*}
X_{t}=M_{t}^{\left( 1\right) }+M_{t}^{\left( 2\right) }+V_{t}
\end{equation*}
where $( M_{t}^{\left( 1\right) }) $ is a continuous local
martingale with respect to a filtration $\left\{ \mathcal{F}_{t}\right\} $, $%
( \hat{M}_{t}^{\left( 2\right) }) $, defined as $\hat{M}%
_{t}^{\left( 2\right) }=M_{T-t}^{\left( 2\right) }$, is a continuous local
martingale with respect to a filtration $\left\{ \mathcal{H}_{t}\right\} $, $%
\left( V_{t}\right) $ is a bounded variation process, $\left( X_{t}\right) $
is adapted to $\left\{ \mathcal{F}_{t}\right\} $ and $( \hat{X}%
_{t}) $ is adapted to $\left\{ \mathcal{H}_{t}\right\} $, and finally
we have 
\begin{equation}
[ M^{\left( 1\right) \alpha }-M^{\left( 2\right) \alpha }] \equiv
0\quad \text{for all }\alpha =1,...,d.  \label{cancell}
\end{equation}
These processes arise in the theory of Dirichlet forms and relevant examples
are the reversible semimartingales. For these processes it is possible to
define stochastic integrals in the sense of Stratonovich, taking advantage
of cancellations coming from assumption (\ref{cancell}). 
For all continuous $1$-forms $\varphi$ we have 
\begin{equation}
\label{lyons-zheng} 
\begin{split}
\int_{0}^{T}\langle \varphi \left( X_{t}\right) ,\circ
dX_{t}\rangle  
 & = \int_{0}^{T}\langle \varphi \left( X_{t}\right) ,dM_{t}^{\left(
1\right) }\rangle -\int_{0}^{T}\langle \varphi ( \hat{X}%
_{t}) ,d\hat{M}_{t}^{\left( 2\right) }\rangle
\\ & \qquad 
+\int_{0}^{T}\left\langle \varphi \left( X_{t}\right) ,dV_{t}\right\rangle     
\end{split}
\end{equation}
where the first two integrals on the right-hand-side are usual It\^{o}
integrals. Indeed, arguing a little bit formally (one has to repeat the
computations on finite sums and for a regularized $\varphi $ and \ prove the
final result taking the limit in probability), we have 
\begin{equation*}
  \begin{split}
    \int_{0}^{T}& \left\langle \varphi \left( X_{t}\right) ,\circ
dX_{t}\right\rangle \\
&= \int_{0}^{T} \langle \varphi \left( X_{t}\right) ,\circ
dM_{t}^{\left( 1\right) } \rangle +\int_{0}^{T} \langle \varphi
\left( X_{t}\right) ,\circ dM_{t}^{\left( 2\right) } \rangle
+\int_{0}^{T} \langle \varphi \left( X_{t}\right) ,
dV_{t} \rangle \\
&= \int_{0}^{T} \langle \varphi \left( X_{t}\right) ,\circ
dM_{t}^{\left( 1\right) } \rangle -\int_{0}^{T}\langle \varphi
( \hat{X}_{t}) ,\circ d\hat{M}_{t}^{\left( 2\right)
}\rangle +\int_{0}^{T}\langle \varphi \left( X_{t}\right)
,dV_{t}\rangle .  
  \end{split}
\end{equation*}
Moreover 
\begin{equation*}
\int_{0}^{T}\langle \varphi \left( X_{t}\right) ,\circ dM_{t}^{\left(
1\right) }\rangle =\sum_{\alpha =1}^{d}\left( \int_{0}^{T}\varphi
_{\alpha }\left( X_{t}\right) dM_{t}^{\left( 1\right) \alpha }+\frac{1}{2}
[ \varphi _{\alpha }\left( X\right) ,M^{\left(
1\right) \alpha }] _{T}\right)
\end{equation*}
\begin{equation*}
\int_{0}^{T}\langle \varphi ( \hat{X}_{t}) ,\circ d\hat{M}%
_{t}^{\left( 2\right) }\rangle =\sum_{\alpha =1}^{d}\left(
\int_{0}^{T}\varphi _{\alpha }( \hat{X}_{t}) d\hat{M}_{t}^{\left(
2\right) \alpha }+\frac{1}{2} [ \varphi _{\alpha }( \hat{X}) ,%
\hat{M}^{\left( 2\right) \alpha }] _{T}\right)
\end{equation*}
\begin{equation*}
\left[ \varphi _{\alpha }\left( \hat{X}\right) ,\hat{M}^{\left( 2\right)
\alpha }\right] _{T}=\left[ \varphi _{\alpha }\left( X\right) ,M^{\left(
2\right) \alpha }\right] _{T}
\end{equation*}
and finally 
\begin{eqnarray*}
[ \varphi _{\alpha }\left( X\right) ,M^{\left( 1\right) \alpha }%
] _{T}-[ \varphi _{\alpha }( \hat{X}) ,\hat{M}^{\left(
2\right) \alpha }] _{T} 
=[ \varphi _{\alpha }\left( X\right) ,M^{\left( 1\right) \alpha
}-M^{\left( 2\right) \alpha }] _{T}=0
\end{eqnarray*}
by (\ref{cancell}). This proves (\ref{lyons-zheng}).

\begin{theorem}
Let $\left( X_{t}\right) $ be a continuous Lyons-Zheng process of the form 
\begin{equation*}
X_{t}=M_{t}^{\left( 1\right) }+M_{t}^{\left( 2\right) }+V_{t}
\end{equation*}
as above. Then the stochastic $1-$current $\varphi \mapsto S\left( \varphi
\right) $ defined by (\ref{strat}) has a pathwise realization $\mathcal{S}%
\left( \omega \right) $, with 
\begin{equation*}
\mathcal{S}\left( \omega \right) \in H^{-s}\left( \mathbb{R}^{d},\mathbb{R}%
^{d}\right) \text{ \qquad }P\text{-a.s.}
\end{equation*}
for all $s>\frac{d}{2}$. If in addition 
\begin{equation*}
\left[ M^{\left( 1\right) }\right] _{T}\;,\;\left[ M^{\left( 2\right) }%
\right] _{T}\;,\;\left\| V\right\| _{var}^{2}\in L^{1}\left( \Omega \right) ,
\end{equation*}
then we also have the integrability property 
\begin{equation*}
\mathcal{S}\left( .\right) \in L^{2}\left( \Omega ,H^{-s}\left( \mathbb{R}%
^{d},\mathbb{R}^{d}\right) \right) .
\end{equation*}
In particular these results hold true for reversible semimartingales.
\bigskip
\end{theorem}

\begin{proof}
Because of (\ref{lyons-zheng}), it is sufficient to prove the result for
each one of the three addenda separately. For the last one it is true by
ordinary integral calculus (recall that $H^{s}\left( \mathbb{R}^{d},\mathbb{R%
}^{d}\right) $ is continuously embedded into the space of continuous $1$%
-forms), while for the first two it is a consequence of the last claim of
the previous theorem on semimartingales. The proof is complete.
\end{proof}

\section{$\protect\gamma $-H\"{o}lder curves ($\protect\gamma >1/2$) define integral flat chains}

The results of this section do not require or involve any stochastic
structure of $\left( X_{t}\right) _{t\in \left[ 0,T\right] }$, which
therefore is supposed to be a deterministic function. The results are of
course applicable to stochastic processes, just path by path (it is also
easy to check that all the quantities constructed below depend measurably on
the random parameter).

Our approach gives a fourmula for the stochastic integral in terms of
double random integrals that seems to be new and could be used in
contexts like the stochastic analisys of fractional Brownian motion.

Let $\left( X_{t}\right) _{t\in \left[ 0,T\right] }$ be a $\gamma $%
-H\"{o}lder continuous function with values in $\mathbb{R}^{d}$, with $\gamma
>1/2$. We first analyse its mollifications. Let $\eta $ be a
real-valued piece-wise $C^{1}$ function on $\mathbb{R}$, with compact
support in $B_{1}\left( 0\right) $, non-negative, such that $\int \eta
\left( t\right) dt=1$. For all $\alpha \in (0,1]$ we set 
\begin{equation*}
\eta _{\alpha }\left( t\right) =\frac{1}{\alpha }\eta \left( \frac{t}{\alpha 
}\right)
\end{equation*}
so that $\int \eta _{\alpha }\left( t\right) dt=1$. We also set 
\begin{equation*}
A=\left[ 0,T\right] \times (0,1]
\end{equation*}
\begin{equation*}
X_{t,\alpha }=( \eta _{\alpha }\ast \tilde{X}) _{t}=\int \eta
_{\alpha }\left( t-s\right) \tilde{X}_{s}\,ds,\qquad \left( t,\alpha \right)
\in A
\end{equation*}
where $( \tilde{X}_{t} ) _{t\in \left[ -1,T+1\right] }$ is any $%
C^{\gamma }$ extension of $X_{t}$. Therefore $\left( X_{t,\alpha }\right)
_{t\in A}$ is a $C^{1}$ mapping from $A$ to $\mathbb{R}^{d}$ (notice that $%
\alpha =0$ is excluded). We shall write $\left( X_{t}\right) $ and $\left(
X_{t,\alpha }\right) $ for shortness.

We have the following formulae and estimates: 
\begin{eqnarray*}
\frac{\partial X_{t,\alpha }}{\partial t} &=&\int \frac{\partial \eta
_{\alpha }}{\partial t}\left( t-s\right) \tilde{X}_{s}\,ds=\left( \frac{%
\partial \eta _{\alpha }}{\partial t}\ast \tilde{X}\right) _{t} \\
&=&\int \frac{\partial \eta _{\alpha }}{\partial t}\left( t-s\right) ( 
\tilde{X}_{s}-\tilde{X}_{t}) \,ds
\end{eqnarray*}
\begin{eqnarray*}
\frac{\partial X_{t,\alpha }}{\partial \alpha } &=&\int \frac{\partial \eta
_{\alpha }}{\partial \alpha }\left( t-s\right) \tilde{X}_{s}\,ds=\left( 
\frac{\partial \eta _{\alpha }}{\partial \alpha }\ast \tilde{X}\right) _{t}
\\
&=&\int \frac{\partial \eta _{\alpha }}{\partial \alpha }\left( t-s\right)
( \tilde{X}_{s}-\tilde{X}_{t}) \,ds
\end{eqnarray*}
\begin{equation*}
\frac{\partial \eta _{\alpha }}{\partial t}\left( t\right) =-\frac{1}{\alpha
^{2}}\eta ^{\prime }\left( \frac{t}{\alpha }\right) ,\quad \frac{\partial
\eta _{\alpha }}{\partial \alpha }\left( t\right) =-\frac{1}{\alpha ^{2}}%
\eta \left( \frac{t}{\alpha }\right) -\frac{t}{\alpha ^{3}}\eta ^{\prime
}\left( \frac{t}{\alpha }\right)
\end{equation*}
\begin{equation*}
\left| \frac{\partial X_{t,\alpha }}{\partial t}\right| \leq C\alpha
^{\gamma -1},\quad \left| \frac{\partial X_{t,\alpha }}{\partial \alpha }%
\right| \leq C\alpha ^{\gamma -1}
\end{equation*}
because 
\begin{equation*}
\begin{split}
\left| \frac{\partial X_{t,\alpha }}{\partial t}\right| &=\left| \int \frac{%
\partial \eta _{\alpha }}{\partial t}\left( t-s\right) ( \tilde{X}_{s}-%
\tilde{X}_{t}) \,ds\right| \\
&\leq \frac{C}{\alpha ^{2}}\int \left| \eta ^{\prime }\left( \frac{t-s}{%
\alpha }\right) \right| \left| t-s\right| ^{\gamma }\,ds \\
&= C\alpha ^{\gamma -1}\int \left| \eta
^{\prime }\left( r\right) \right| \left| r\right| ^{\gamma }\,dr
\qquad \text{with $r=\frac{t-s}{\alpha }$}
\end{split}
\end{equation*}
\begin{eqnarray*}
\left| \frac{\partial X_{t,\alpha }}{\partial \alpha }\right| &=&\left| \int 
\frac{\partial \eta _{\alpha }}{\partial \alpha }\left( t-s\right) \left( 
\tilde{X}_{s}-\tilde{X}_{t}\right) \,ds\right| \\
&\leq &\frac{C}{\alpha ^{2}}\left| \int \left| \eta \left( \frac{t-s}{\alpha 
}\right) +\frac{t-s}{\alpha }\eta ^{\prime }\left( \frac{t-s}{\alpha }%
\right) \right| \left| t-s\right| ^{\gamma }\,ds\right| \\
&=& C\alpha ^{\gamma -1}\left| \int \left|
\eta \left( r\right) +r\eta ^{\prime }\left( r\right) \right| \left|
r\right| ^{\gamma }\,dr\right| \qquad \text{with $r=\frac{t-s}{\alpha }$}.
\end{eqnarray*}
The next theorem is perhaps the main result of this paper. It states that
any $\gamma $-H\"{o}lder continuous curve in $\mathbb{R}^{d}$, with $\gamma >%
1/2$, defines an integral flat chain, denoted by $%
\int_{0}^{T}\left\langle \varphi \left( X_{t}\right) ,\,\circ
dX_{t}\right\rangle $ and gives a formula in terms of a double
integral. 
For curves with such regularity there is no distinction between Stratonovich and It\^{o} integrals (the quadratic
variation is easily proved to be zero), so we keep the notation of
Stratonovich integral since we believe it should be the appropriate one in
case the following theorem will have some kind of generalization to
semimartingales in the future. Notice that the following result is not
entirely trivial a priori since $\left( X_{t}\right) $ may have infinite
mass and Hausdorff dimension $\gamma ^{-1}$.

\begin{theorem}
\label{th:prop_key1}
Assume that $\left( X_{t}\right) $ is a $\gamma $-H\"{o}lder continuos curve
in $\mathbb{R}^{d}$, with $\gamma >1/2$. Then:

a) $\left( X_{t,\alpha }\right) $ defines the following i.m. rectifiable $2$%
-current $T$: 
\begin{eqnarray*}
T\left( \psi \right)  &=&\int_{A}\left\langle \psi \left( X_{t,\alpha
}\right) ,\frac{\partial X_{t,\alpha }}{\partial t}\wedge \frac{\partial
X_{t,\alpha }}{\partial \alpha }\right\rangle \,dtd\alpha  \\
&=&\int_{A}\sum_{i<j}\psi _{ij}\left( X_{t,\alpha }\right) \left( \frac{%
\partial X_{t,\alpha }^{j}}{\partial t}\frac{\partial X_{t,\alpha }^{i}}{%
\partial \alpha }-\frac{\partial X_{t,\alpha }^{i}}{\partial t}\frac{%
\partial X_{t,\alpha }^{j}}{\partial \alpha }\right) \,dtd\alpha 
\end{eqnarray*}
for all continuous $2$-forms $\psi $ on $\mathbb{R}^{d}$ represented as 
\begin{equation*}
\sum_{i<j}\psi _{ij}\left( x\right) \,dx^{i}\wedge dx^{j}.
\end{equation*}

b) for every continuously differentiable $1$-form $\varphi $, the following
limit 
\begin{equation*}
\lim_{\alpha \rightarrow 0}\int_{0}^{T}\left\langle \varphi \left(
X_{t,\alpha }\right) ,\circ \,dX_{t,\alpha }\right\rangle 
\end{equation*}
exists (the integral is understood as a classical Riemann integral) and will
be denoted by $\int_{0}^{T}\left\langle \varphi \left( X_{t}\right) ,\circ
\,dX_{t}\right\rangle $.

c) The mapping $\varphi \mapsto \int_{0}^{T}\left\langle \varphi \left(
X_{t}\right) ,\circ \,dX_{t}\right\rangle $ is the integral flat chain (of
degree one) given by
\begin{equation}
  \label{chain}
  \begin{split}
\int_{0}^{T} &\left\langle \varphi \left( X_{t}\right) ,\circ
\,dX_{t}\right\rangle  = T\left( d\varphi \right)   
+\int_{0}^{T}\left\langle \varphi \left( X_{t,1}\right) ,\,\frac{\partial
X_{t,1}}{\partial t}\right\rangle \,dt   
\\ & \quad -\int_{0}^{1}\left\langle \varphi \left( X_{T,\alpha }\right) ,\,\frac{%
\partial X_{T,\alpha }}{\partial \alpha }\right\rangle \,d\alpha
+\int_{0}^{1}\left\langle \varphi \left( X_{0,\alpha }\right) ,\,\frac{%
\partial X_{0,\alpha }}{\partial \alpha }\right\rangle \,d\alpha     
  \end{split}
\end{equation}
for all continuously differentiable $1$-forms $\varphi $ on $\mathbb{R}^{d}$
represented as 
$
\sum_{i=1}^{d}\varphi _{i}\left( y\right) \,dy^{i}
$ and where
\begin{equation*}
T\left( d\varphi \right) =\int_{A}\sum_{i,j=1}^{d}\frac{\partial \varphi _{i}%
}{\partial y_{j}}\left( X_{t,\alpha }\right) \left( \frac{\partial
X_{t,\alpha }^{j}}{\partial t}\frac{\partial X_{t,\alpha }^{i}}{\partial
\alpha }-\frac{\partial X_{t,\alpha }^{i}}{\partial t}\frac{\partial
X_{t,\alpha }^{j}}{\partial \alpha }\right) \,dtd\alpha .
\end{equation*}
\end{theorem}

\begin{proof}
\textbf{Step 1} (The support $\mathcal{M}$ of $T$ is a $2$-rectifiable set).
Let us introduce the function $f:A\rightarrow \mathbb{R}^{d}$ defined as 
$
f(t,\alpha )=X_{t,\alpha }
$
and let $\mathcal{M}=f(A)$. We have 
\begin{equation*}
Df=\left( 
\begin{array}{cc}
\displaystyle \frac{\partial X_{t,\alpha }}{\partial t} & \displaystyle\frac{\partial X_{t,\alpha }}{%
\partial \alpha }
\end{array}
\right) 
\end{equation*}
and 
\begin{equation*}
\left( Df\right) ^{\ast }Df=\left( 
\begin{array}{cc}
\displaystyle\left| \frac{\partial X_{t,\alpha }}{\partial t}\right| ^{2} & \displaystyle\left\langle 
\frac{\partial X_{t,\alpha }}{\partial t},\frac{\partial X_{t,\alpha }}{%
\partial \alpha }\right\rangle  \\ 
\displaystyle\left\langle \frac{\partial X_{t,\alpha }}{\partial t},\frac{\partial
X_{t,\alpha }}{\partial \alpha }\right\rangle  & \displaystyle\left| \frac{\partial
X_{t,\alpha }}{\partial \alpha }\right| ^{2}
\end{array}
\right) 
\end{equation*}
\begin{eqnarray*}
J_{f} &=&\sqrt{\det \left( Df\right) ^{\ast }Df} \\
&=&\sqrt{\left| \frac{\partial X_{t,\alpha }}{\partial \alpha }\right|
^{2}\left| \frac{\partial X_{t,\alpha }}{\partial t}\right|
^{2}-\left\langle \frac{\partial X_{t,\alpha }}{\partial t},\frac{\partial
X_{t,\alpha }}{\partial \alpha }\right\rangle ^{2}}.
\end{eqnarray*}
By the estimates of the previous subsection we have 
\begin{equation*}
J_{f}\leq C\alpha ^{2\left( \gamma -1\right) }
\end{equation*}
which is integrable on $A$: 
\begin{equation*}
\int_{A}J_{f}\,dtd\alpha <\infty .
\end{equation*}
This is the basic property that will imply the final result.

Let $A_{n}$ be defined as 
\begin{equation*}
A_{n}=\left\{ \left( t,\alpha \right) \in \mathbb{R}^{2}:t\in \left[ 0,T%
\right] ,\;\alpha \in \left( \frac{1}{n+1},1\right] \right\} .
\end{equation*}
We have $\mathcal{M}=\,\bigcup_{n}f\left( A_{n}\right) $ where $f\lfloor
A_{n}$ are Lipschitz functions. In particular (see \cite{GMS} p. 75), $%
f\left( A_{n}\right) $ are $\mathcal{H}^{2}$-measurable sets, and so it is $%
\mathcal{M}$. Since $f\lfloor A_{n}$ is Lipschitz, the area formula gives us 
\begin{equation}
\int_{f\left( A_{n}\right) }\mathcal{H}^{0}\left( f^{-1}\left( x\right) \cap
A_{n}\right) d\mathcal{H}^{2}\left( x\right) =\int_{A_{n}}J_{f}\,dtd\alpha
\leq \int_{A}J_{f}\,dtd\alpha <\infty .  \label{areaformula}
\end{equation}
On one side it follows that 
\begin{equation*}
\mathcal{H}^{2}\left( f\left( A_{n}\right) \right) =\int_{f\left(
A_{n}\right) }d\mathcal{H}^{2}\left( x\right) \leq \int_{A}J_{f}\,dtd\alpha
<\infty 
\end{equation*}
for all $n$, and therefore, by the continuity of $\mathcal{H}^{2}$, 
\begin{equation*}
\mathcal{H}^{2}\left( \mathcal{M}\right) <\infty .
\end{equation*}
This completes the proof that $\mathcal{M}$ is a $2$-rectifiable set in $%
\mathbb{R}^{d}$.

\textbf{Step 2} ($T$ is an i.m. rectifiable $2$-current). Again from (\ref
{areaformula}) we deduce that the function $N\left( f,A,.\right) :\mathcal{M}%
\rightarrow \mathbb{R}$ defined as $N\left( f,A,x\right) =\mathcal{H}%
^{0}\left( f^{-1}\left( x\right) \right) $ is $\mathcal{H}^{2}\lfloor 
\mathcal{M}$-a.s. finite and integer valued, and also $\mathcal{H}%
^{2}\lfloor \mathcal{M}$-integrable. Indeed, by the monotone convergence
theorem 
\begin{eqnarray}
\int_{\mathcal{M}}N\left( f,A,x\right) \,d\mathcal{H}^{2}\left( x\right) 
&=&\lim_{n\rightarrow \infty }\int_{f\left( A_{n}\right) }\mathcal{H}%
^{0}\left( f^{-1}\left( x\right) \cap A_{n}\right) d\mathcal{H}^{2}\left(
x\right)   \label{N} \\
&\leq &\int_{A}J_{f}\,\,dtd\alpha <\infty .  \notag
\end{eqnarray}
In particular, it follows that the following vector field is well defined $%
\mathcal{H}^{2}\lfloor \mathcal{M}$-a.s. since it is given by a finite sum: 
$$
v\left( x\right) =\sum_{\left( t,\alpha \right) \in f^{-1}\left( x\right) }%
\frac{\zeta_{t,\alpha}}{|\zeta_{t,\alpha}|}
$$
with
$$
\zeta_{t,\alpha} := 
  \left( \frac{\partial
      X_{t,\alpha }}{\partial t}\wedge \frac{\partial X_{t,\alpha }}{\partial
      \alpha }\right)
.
$$


Let us introduce the $2$-current $T_{n}$ defined as
\begin{equation*}
T_{n}\left( \psi \right) :=\int_{A_{n}}\left\langle \psi \left( X_{t,\alpha
}\right) ,\zeta_{t,\alpha}\right\rangle \,\,dtd\alpha .
\end{equation*}
Since $f\lfloor A_{n}$ is Lipschitz, from the change of variable formula ( 
\cite{GMS} p. 75) we have 
\begin{equation*}
  \begin{split}
T_{n}\left( \psi \right) & =\int_{A_{n}}\left\langle \psi \left( X_{t,\alpha
}\right) ,\frac{\zeta_{t,\alpha}}{|\zeta_{t,\alpha}|}%
\right\rangle J_{f}\,\,dtd\alpha 
\\ &= \int_{f\left( A_{n}\right) }\sum_{\left( t,\alpha \right) \in
f^{-1}\left( x\right) \cap A_{n}}\left\langle \psi \left( X_{t,\alpha
}\right) ,\frac{\zeta_{t,\alpha}}{|\zeta_{t,\alpha}|}
\right\rangle \,\,d\mathcal{H}^{2}\left( x\right)  
\\&=\int_{f\left( A_{n}\right) }\left\langle \psi \left( x\right)
,v_{n}\left( x\right) \right\rangle \,\,d\mathcal{H}^{2}\left( x\right) 
  \end{split}
\end{equation*}
where $v_{n}\left( x\right) $ is defined as
\begin{equation*}
v_{n}\left( x\right) :=\sum_{\left( t,\alpha \right) \in f^{-1}\left(
x\right) \cap A_{n}}\frac{\zeta_{t,\alpha}}{|\zeta_{t,\alpha}|}
\end{equation*}
Since $f(A_{n})$ is a $2$-rectifiable set, the density $\Theta ^{2}\left(
f(A_{n}),x\right) $ is equal to $1$ for $\mathcal{H}^{2}\lfloor f(A_{n})$%
-a.e. $x$ (see \cite{Fed}, th. 3.2.19), hence the cone of approximate
tangent vectors of $f(A_{n})$ at such $x$ is the plane identified by any one
of the $2$-vectors $\zeta_{t,\alpha}/|\zeta_{t,\alpha}|$. These $2$-vectors are unitary, hence they
coincide up to the constant $\pm 1$. Therefore there exist an $\mathcal{H}%
^{2}$-measurable unitary $2$-vector field $\xi _{n}$ on $f(A_{n})$, with $%
\xi _{n}\left( x\right) \in T_{x}f(A_{n})$ for $\mathcal{H}^{2}\lfloor
f(A_{n})$-a.e. $x,$ and an integer valued and $\mathcal{H}^{2}\lfloor
f(A_{n})$ integrable (by (\ref{N})) function $\theta
_{n}:f(A_{n})\rightarrow \mathbb{R}$, such that $v_{n}=\xi _{n}\theta _{n}$ $%
\mathcal{H}^{2}\lfloor f(A_{n})$-a.s. This proves the (well known) fact that 
$T_{n}$ is a i.m. rectifiable $2$-current:
\begin{equation}
T_{n}\left( \psi \right) =\int_{f\left( A_{n}\right) }\left\langle \psi ,\xi
_{n}\right\rangle \theta _{n}\,\,d\mathcal{H}^{2}.  \label{Tn}
\end{equation}

Let us now decompose $\mathcal{M}$ into the sets $\mathcal{N}_{n}=f\left(
A-A_{n}\right) $ and $\mathcal{M}_{n}=f\left( A_{n}\right) \setminus 
\mathcal{N}_{n}$. We have
\begin{eqnarray*}
\mathcal{H}^{2}\left( \mathcal{N}_{n}\right)  &=&\int_{A\setminus A_{n}}d%
\mathcal{H}^{2}\left( x\right)  \\
&\leq &\int_{f\left( A\setminus A_{n}\right) }\mathcal{H}^{0}\left(
f^{-1}\left( x\right) \cap \left( A\setminus A_{n}\right) \right) \,d%
\mathcal{H}^{2}\left( x\right)  \\
&=&\int_{A\setminus A_{n}}J_{f}\,\,dtd\alpha \rightarrow 0\text{ as }%
n\rightarrow \infty \text{.}
\end{eqnarray*}
In a sense, the subset $\mathcal{M}_{n}$ of $\mathcal{M}$ is a stable part
when $n$ increases (it is not visited any more by $f$), and the remainder $%
\mathcal{N}_{n}$ becomes negligible. We have, for all $m>n$,  
\begin{equation}
T_{x}\mathcal{M}=T_{x}f\left( A_{m}\right) =T_{x}f\left( A_{n}\right) \text{
for a.e. }x\in \mathcal{M}_{n}\text{.}  \label{tangent}
\end{equation}
Indeed, arguing again on the density (\cite{Fed}, th. 3.2.19), for a.e. $%
x\in \mathcal{M}_{n}$ the three tangent planes (that exist) are the only
part of the corresponding cones of approximate tangent vectors (a priori one
bigger than the other), so they coincide. We also have (directly from the
definitions)
\begin{equation*}
v_{m}\left( x\right) =v_{n}\left( x\right) \text{ for a.e. }x\in \mathcal{M}%
_{n}.
\end{equation*}
Choosing $\theta _{n}$ above with non negative values (it is always possible),
we deduce
\begin{equation*}
\xi _{m}\left( x\right) =\xi _{n}\left( x\right) \text{ and }\theta
_{m}\left( x\right) =\theta _{n}\left( x\right) \ \text{for a.e. }x\in 
\mathcal{M}_{n}.
\end{equation*}
It is therefore well defined an $\mathcal{H}^{2}$-measurable unitary $2$%
-vector field $\xi $ on $\mathcal{M}$, with (by (\ref{tangent}))  $\xi
\left( x\right) \in T_{x}\mathcal{M}$ for $\mathcal{H}^{2}\lfloor \mathcal{M}
$-a.e. $x$, and an integer valued and $\mathcal{H}^{2}\lfloor \mathcal{M}$
integrable (again by (\ref{N})) function $\theta :\mathcal{M}\rightarrow 
\mathbb{R}$, such that
\begin{equation*}
\xi \left( x\right) =\xi _{n}\left( x\right) \text{ and }\theta \left(
x\right) =\theta _{n}\left( x\right) \ \text{for a.e. }x\in \mathcal{M}_{n}.
\end{equation*}
From (\ref{Tn}) we have
\begin{equation*}
T_{n}\left( \psi \right) =\int_{\mathcal{M}_{n}}\left\langle \psi ,\xi
\right\rangle \theta \,\,d\mathcal{H}^{2}+\int_{\mathcal{N}_{n}}\left\langle
\psi ,v_{n}\right\rangle \,\,d\mathcal{H}^{2}.
\end{equation*}
Since $\int_{A}J_{f}\,\,dtd\alpha <\infty $, $T_{n}\left( \psi \right)
\rightarrow T\left( \psi \right) $ as $n\rightarrow \infty $. By monotone
convergence, the second term converges to $\int_{\mathcal{M}}\left\langle
\psi ,\xi \right\rangle \theta \,\,d\mathcal{H}^{2}$. Finally,
\begin{eqnarray*}
\left| \int_{\mathcal{N}_{n}}\left\langle \psi ,v_{n}\right\rangle \,\,d%
\mathcal{H}^{2}\right|  &\leq &\int_{\mathcal{N}_{n}}\left\| \psi \right\|
_{\infty }\left\| v\left( x\right) \right\| \,\,d\mathcal{H}^{2}\left(
x\right)  \\
&\leq &\left\| \psi \right\| _{\infty }\int_{\mathcal{N}_{n}}N\left(
f,A,x\right) \,\,d\mathcal{H}^{2}\left( x\right) 
\end{eqnarray*}
that converges to zero by (\ref{N}). Therefore in the limit we obtain
\begin{equation*}
T\left( \psi \right) =\int_{\mathcal{M}}\left\langle \psi ,\xi \right\rangle
\theta \,\,d\mathcal{H}^{2}
\end{equation*}
so $T$ is an i.m. rectifiable $2$-current.

\textbf{Step 3} (conclusion). The proof of b) and of formula (\ref{chain})
is elementary: first one writes a formula like (\ref{chain}) for classical
integrals on $A_{\alpha }=\left[ 0,T\right] \times (\alpha ,1]$ (the
classical Stokes formula), then pass to the limit. The last three integrals
(boundary terms) in (\ref{chain}) are well defined classical integrals and
define i.m. rectifiable $1$-currents. Therefore (\ref{chain}) tells us that $%
\int_{0}^{T}\left\langle \varphi \left( X_{t}\right) ,\circ
\,dX_{t}\right\rangle $ is an integral flat chain (recall that $T\left(
d\varphi \right) =\partial T\left( \varphi \right) $). The proof is complete.
\end{proof}

A priori the integral flat chain $\int_{0}^{T}\left\langle \varphi \left(
X_{t}\right) ,\,\circ dX_{t}\right\rangle $ defined above may depend
on the choice of $%
\eta $ and of the extension $\tilde{X}$. This is not the case, as the
following proposition asserts.

\begin{proposition}
\label{th:prop_key2}
a) If $( X_{t}^{\left( n\right) } ) $ is a sequence of functions
in $C^{\gamma }$, with $\gamma >1/2$, that converges to $\left(
X_{t}\right) \in $ $C^{\gamma }$ in the $C^{\gamma }$ seminorm and in $C^{0}$%
, then for every continuously differentiable $1$-form $\varphi $ we have 
\begin{equation*}
\int_{0}^{T}\langle \varphi ( X_{t}^{\left( n\right) } )
, \circ dX_{t}^{\left( n\right) }\rangle \rightarrow
\int_{0}^{T}\left\langle \varphi \left( X_{t}\right) ,\,\circ
dX_{t}\right\rangle 
\end{equation*}
where the $C^{\gamma }$- extensions of $( X_{t}^{\left( n\right)
}) $ and $\left( X_{t}\right) $ are arbitrarily chosen and the same $%
\eta $ is taken.

b) If $\left( X_{t}\right) $ is continuously differentiable, $%
\int_{0}^{T}\left\langle \varphi \left( X_{t}\right) ,\,\circ
dX_{t}\right\rangle $ defined above coincides with the usual Riemann
integral, for all $\eta $ and $C^{\gamma }$- extension of $\left(
X_{t}\right) $.

c) If $\left( X_{t}\right) $ is in $C^{\gamma }$ with $\gamma >1/2$%
, the definition of the integral flat chain $\int_{0}^{T}\left\langle
\varphi \left( X_{t}\right) ,\,\circ dX_{t}\right\rangle $ is independent of 
the choice of $\eta $ and of the $C^{\gamma }$- extension of $\left( X_{t}\right) $.
\end{proposition}
\begin{proof}
a) This follows from the explicit representation formula for the
integral $%
\int_{0}^{T}\langle \varphi( X_{t}^{\left( n\right) })
,\,\circ dX_{t}^{\left( n\right) }\rangle $ and $\int_{0}^{T}\left%
\langle \varphi \left( X_{t}\right) ,\,\circ dX_{t}\right\rangle $, after
some tedious but elementary estimates similar to those used above.

b) It is just Stokes formula.

c) It follows from a) and b) by taking as $( X_{t}^{\left( n\right)
}) $ a sequence of $C^{1}$ approximations of $\left( X_{t}\right) $.
\end{proof}

\begin{remark}
For processes with strong fluctuations, like the fractional Brownian motion,
the following fact is plausible: in $d=2$ the set $f(A)$ contains a lot of
superpositions; in $d\geq 3$ there are only selfintersections in $f(A)$, and
their set has $2$-dimensional Hausdorff measure (in other words, $N(f,A,x)=1$%
, $\mathcal{H}^{2}\lfloor \mathcal{M}$-a.s.).
\end{remark}

\begin{remark}
\label{rem:sobolev}
Recall the spaces $W^{s,p}$ as defined for instance in \cite{Ad}. 
Let 
$W_{C}^{1/2,2}$ be the set of all continuous functions $%
\left( X_{t}\right) _{t\in \left[ 0,T\right] }$ taking values in $\mathbb{R}%
^{d}$ such that 
\begin{equation*}
\int_{0}^{T}\int_{0}^{T}\frac{\left| X_{t}-X_{s}\right| ^{2}}{\left|
t-s\right| ^{2}}dtds+\int_{0}^{T}\frac{\left| X_{T}-X_{s}\right| }{\left|
T-s\right| }ds+\int_{0}^{T}\frac{\left| X_{0}-X_{s}\right| }{\left| s\right| 
}ds<\infty .
\end{equation*}
We easily have 
\begin{equation*}
\bigcup_{\gamma >1/2}C^{\gamma }\left( \left[ 0,T\right] \right)
\subset W_{C}^{1/2,2}
\end{equation*}
The results of the previous section extends to the class
 $W^{1/2,2}_C$:  one can prove that


\begin{equation*}
\int_{0}^{T}\int_{0}^{1}\left( \left| \frac{\partial X_{t,\alpha }}{\partial
t}\right| ^{2}+\left| \frac{\partial X_{t,\alpha }}{\partial \alpha }\right|
^{2}\right) dtd\alpha <\infty
\end{equation*}
\begin{equation*}
\int_{0}^{1}\left( \left| \frac{\partial X_{0,\alpha }}{\partial \alpha }%
\right| +\left| \frac{\partial X_{T,\alpha }}{\partial \alpha }\right|
\right) d\alpha <\infty ;
\end{equation*}
and therefore
\begin{equation*}
\int_{A}J_{f}\,dtd\alpha <\infty .
\end{equation*}


The extension to $W^{1/2,2}_C$ of Theorem~\ref{th:prop_key1} can be
interpreted as a result of trace theory in Sobolev spaces.  
We think that the constructive proof we outline
is simple and straighforward enough to dispense with abstract
arguments, which in addition would require a careful investigation of
the boundary terms at $t=0,T$. 
Moreover as it will be shown in~\cite{GubWIP} the explicit
computations we make give some insights on the difficulties of
extending the result to more irregular paths (e.g. samples of
Brownian motion).
\end{remark}

\begin{remark}
Given two curves $\left( X_{t}\right) $ and $\left( Y_{t}\right) $ of class 
$C^{\gamma }$, $\gamma >1/2$ (the same can be said for the class $%
W_{C}^{1/2,2}$, see remark~\ref{rem:sobolev} above), we may now define the \emph{integral
flat distance} $d_{\mathcal{F}}\left( X,Y\right) $ between them as the
infimum of the numbers 
\begin{equation*}
\mathbf{M}\left( T\right) +\mathbf{M}\left( R\right) \,
\end{equation*}
over all i.m. rectifiable $2$-currents $T$ and i.m. rectifiable $1$-currents 
$R$ such that 
\begin{equation*}
\int_{0}^{T}\left\langle \varphi \left( X_{t}\right) ,\,\circ
dX_{t}\right\rangle -\int_{0}^{T}\left\langle \varphi \left( Y_{t}\right)
,\,\circ dY_{t}\right\rangle =\partial T+R
\end{equation*}
(by the previous theorem, at least one pair $\left( T,R\right) $ exists).
Here $\mathbf{M}\left( .\right) $ denotes the mass of the corresponding
current (see \cite{GMS}). In a sense, this distance is more geometric than
the usual ones employed in stochastic analysis, it is independent of the
parametrization, and thus we hope it may have applications in future
researches.
\end{remark}

\subsection{A $d+1$ dimensional variant}

Geometrically it may appear unpleasant that the support of the rectifiable
current $\left( X_{t,\alpha }\right) $ may have very complex overlaps
(expecially in dimension 2). There are many ways to associate to $\left(
X_{t}\right) $ a piece of boundary of a rectifiable 2-current, some of them
giving us surfaces with a nicer appearance. An easy way is to see $\left(
X_{t}\right) $ in $\mathbb{R}^{d+1}$, namely to consider the function 
\begin{equation*}
t\mapsto Y_{t,\alpha }:=\left( t,X_{t}\right) \text{.}
\end{equation*}
We conjecture that in some relevant examples this mapping is injective.

This modification has the additional advantage that it defines integrals
over 1-forms depending also on $t$. In other words, we define integrals of
the form 
\begin{equation*}
\int_{0}^{T}\left\langle \varphi \left( t,X_{t}\right) \,,\circ
dX_{t}\right\rangle .
\end{equation*}

\begin{theorem}
Assume that $\left( X_{t}\right) $ is a $\gamma $-H\"{o}lder continuos curve
in $\mathbb{R}^{d}$ with $\gamma >1/2$. Then:

a) $\left( Y_{t,\alpha }\right) $ defines the following integer multiplicity
rectifiable $2$-current $\tilde{T}$ in $\mathbb{R}^{d+1}$: 
\begin{eqnarray*}
\tilde{T}\left( \psi \right)  &=&\int_{A}\left\langle \psi \left(
t,X_{t,\alpha }\right) ,\left( 1,\frac{\partial X_{t,\alpha }}{\partial t}%
\right) \wedge \left( 0,\frac{\partial X_{t,\alpha }}{\partial \alpha }%
\right) \right\rangle \,dtd\alpha  \\
&=&\int_{A}\sum_{i,j=1}^{d}\psi _{ij}\left( t,X_{t,\alpha }\right) \left( 
\frac{\partial X_{t,\alpha }^{j}}{\partial t}\frac{\partial X_{t,\alpha }^{i}%
}{\partial \alpha }-\frac{\partial X_{t,\alpha }^{i}}{\partial t}\frac{%
\partial X_{t,\alpha }^{j}}{\partial \alpha }\right) \,dtd\alpha  \\
&&-\int_{A}\sum_{j=1}^{d}\psi _{0j}\left( t,X_{t,\alpha }\right) \frac{%
\partial X_{t,\alpha }^{j}}{\partial \alpha }\,dtd\alpha
+\int_{A}\sum_{i=1}^{d}\psi _{i0}\left( t,X_{t,\alpha }\right) \frac{%
\partial X_{t,\alpha }^{i}}{\partial \alpha }\,dtd\alpha 
\end{eqnarray*}
for all continuous $2$-forms $\psi $ on $\mathbb{R}^{d+1}$ represented as 
\begin{equation*}
\sum_{i,j=1}^{d}\psi _{ij}\left( t,y\right) \,dy^{j}\wedge
dy^{i}+\sum_{j=1}^{d}\psi _{0j}\left( t,y\right) \,dy^{j}\wedge
dt+\sum_{i=1}^{d}\psi _{i0}\left( t,y\right) \,dt\wedge dy^{i}.
\end{equation*}

b) for every continuously differentiable $1$-form $\varphi $, the following
limit 
\begin{equation*}
\lim_{\alpha \rightarrow 0}\int_{0}^{T}\left\langle \varphi \left(
t,X_{t,\alpha }\right) ,\circ \,dX_{t,\alpha }\right\rangle 
\end{equation*}
exists (the integral is understood as a classical Riemann integral) and will
be denoted by $\int_{0}^{T}\left\langle \varphi \left( t,X_{t}\right) ,\circ
\,dX_{t}\right\rangle $.

c) The mapping $\varphi \mapsto \int_{0}^{T}\left\langle \varphi \left(
t,X_{t}\right) ,\circ \,dX_{t}\right\rangle $ is the integral flat chain (of
degree one) given by
\begin{equation*}
  \begin{split}
\int_{0}^{T} &\left\langle \varphi \left( t,X_{t}\right) \,,\circ
dX_{t}\right\rangle  =\tilde{T}\left( d\varphi \right) 
+\int_{0}^{T}\left\langle \varphi \left( t,X_{t,1}\right) ,\,\left( 1,\frac{%
\partial X_{t,1}}{\partial t}\right) \right\rangle \,dt \\
& \quad -\int_{0}^{1}\left\langle \varphi \left( t,X_{T,\alpha }\right) ,\,\left( 0,%
\frac{\partial X_{T,\alpha }}{\partial \alpha }\right) \right\rangle
\,d\alpha +\int_{0}^{1}\left\langle \varphi \left( t,X_{0,\alpha }\right)
,\,\left( 0,\frac{\partial X_{0,\alpha }}{\partial \alpha }\right)
\right\rangle \,d\alpha    
  \end{split}
\end{equation*}
for all $1$-forms $\varphi $ on $\mathbb{R}^{d}$ represented as 
$
\sum_{i=1}^{d}\varphi _{i}\left( t,y\right) \,dy^{i}+\varphi _{0}\left(
t,y\right) \,dt.
$
\end{theorem}

The proof is similar to that of Theorem~\ref{th:prop_key1} in the
previous subsection.  In this case we
define 
\begin{equation*}
\tilde{f}(t,\alpha )=\left( t,X_{t,\alpha }\right) =\left( t,\left( \eta
_{\alpha }\ast X\right) _{t}\right)
\end{equation*}
and we use the formulae 
\begin{equation*}
D\tilde{f}=\left( 
\begin{array}{cc}
\displaystyle 1 & \displaystyle 0 \\ 
\displaystyle\frac{\partial X_{t,\alpha }}{\partial t} & \displaystyle\frac{\partial X_{t,\alpha }}{%
\partial \alpha }
\end{array}
\right)
\end{equation*}
\begin{equation*}
\left( D\tilde{f}\right) ^{\ast }D\tilde{f}=\left( 
\begin{array}{cc}
\displaystyle1+\left| \frac{\partial X_{t,\alpha }}{\partial t}\right| ^{2} &\displaystyle 
\left\langle \frac{\partial X_{t,\alpha }}{\partial t},\frac{\partial
X_{t,\alpha }}{\partial \alpha }\right\rangle \\ 
\displaystyle
\left\langle \frac{\partial X_{t,\alpha }}{\partial t},\frac{\partial
X_{t,\alpha }}{\partial \alpha }\right\rangle & \displaystyle\left| \frac{\partial
X_{t,\alpha }}{\partial \alpha }\right| ^{2}
\end{array}
\right)
\end{equation*}
\begin{eqnarray*}
J_{\tilde{f}} &=&\sqrt{\det \left( D\tilde{f}\right) ^{\ast }D\tilde{f}} \\
&=&\sqrt{\left| \frac{\partial X_{t,\alpha }}{\partial \alpha }\right|
^{2}\left( 1+\left| \frac{\partial X_{t,\alpha }}{\partial t}\right|
^{2}\right) -\left\langle \frac{\partial X_{t,\alpha }}{\partial t},\frac{%
\partial X_{t,\alpha }}{\partial \alpha }\right\rangle ^{2}}.
\end{eqnarray*}
Moreover, the analog of Proposition~\ref{th:prop_key2} holds true.

\subsection{Remarks for Brownian curves}
The results of the previous section do not apply to the typical paths of a
Brownian motion $\left( W_{t}\right) $ in $\mathbb{R}^{d}$
  and presumably
 it is not possible to modify the argument and prove that such paths are
 integral flat currents. 
To clarify this point we report below a negative result on the
summability of  $J_{f}$ in the case of Brownian motion.

\begin{lemma}
We have identities in law (for every $\lambda >0$) between the following
processes, in the range $t\geq 1$, $\alpha \in (0,1]$: 
\begin{equation*}
\left( \frac{\partial W_{t,\alpha }}{\partial t}\left( \lambda t,\lambda
\alpha \right) ,\frac{\partial W_{t,\alpha }}{\partial \alpha }\left(
\lambda t,\lambda \alpha \right) \right) \overset{\mathcal{L}}{=}\frac{1}{%
\sqrt{\lambda }}\left( \frac{\partial W_{t,\alpha }}{\partial t}\left(
t,\alpha \right) ,\frac{\partial W_{t,\alpha }}{\partial \alpha }\left(
t,\alpha \right) \right)
\end{equation*}
\begin{equation*}
J_{f}\left( \lambda t,\lambda \alpha \right) \overset{\mathcal{L}}{=}\frac{1%
}{\lambda }J_{f}\left( t,\alpha \right) .
\end{equation*}
Moreover, the law of these processes do not change by a time shift.
\end{lemma}

\begin{proof}
Recall that the law of the process $\left( W_{\lambda t}\right) $ is equal
to ther law of $( \sqrt{\lambda }W_{t}) $. Then 
\begin{eqnarray*}
\frac{\partial W_{t,\alpha }}{\partial t}\left( \lambda t,\lambda \alpha
\right) &=&\int \frac{\partial \eta _{\alpha }}{\partial t}\left( \lambda
\left( t-s\right) ,\lambda \alpha \right) W_{\lambda s}ds \\
\overset{\mathcal{L}}{=-}\int \frac{1}{\left( \lambda \alpha \right) ^{2}}%
\eta ^{\prime }\left( \frac{t-s}{\alpha }\right) W_{s}\lambda ^{\frac{3}{2}%
}ds &=&\frac{1}{\sqrt{\lambda }}\frac{\partial W_{t,\alpha }}{\partial t}%
\left( t,\alpha \right) .
\end{eqnarray*}
The computation for $\partial W_{t,\alpha }/\partial \alpha$ is
similar. 
To give a complete proof one has to repeat
the previous computation jointly on the two previous random variables and
jointly at generic points $\left( t_{1},\alpha _{1}\right) $, \dots , $\left(
t_{n},\alpha _{n}\right) $, which is just more lengthy. The property for $%
J_{f}$ is a direct consequence of the previous result. The invariance under
time shift is due to the representation (obtained by integration by parts) 
\begin{equation*}
\frac{\partial W_{t,\alpha }}{\partial t}=\int \eta _{\alpha }\left(
t-s\right) \,dW_{s}
\end{equation*}
\begin{equation*}
\frac{\partial W_{t,\alpha }}{\partial \alpha }=-\int \eta _{\alpha }\left(
t-s\right) \frac{t-s}{\alpha }\,dW_{s}
\end{equation*}
along with the stationarity of the increments of the Brownian motion.
\end{proof}

\begin{proposition}
Over the set $(t,\alpha )\in \left[ 1,2\right] \times (0,1]$ we have 
\begin{equation*}
E\int_{\left[ 1,2\right] \times (0,1]}J_{f}\,dtd\alpha =+\infty .
\end{equation*}
\end{proposition}

\begin{proof}
Split $\left[ 1,2\right] \times (0,1]$ in strips $S_{n}=\left[ 1,2\right]
\times (2^{-n-1},2^{-n}]$, $n=0,1,$... Split each strip in
squares: 
\begin{equation*}
S_{n}=\bigcup_{k=0}^{2^{n+1}-1}S_{n,k}\text{ with }S_{n,k}=\left[ 1+\frac{k}{%
2^{n+1}},1+\frac{k+1}{2^{n+1}}\right] \times (\frac{1}{2^{n+1}},\frac{1}{%
2^{n}}].
\end{equation*}
Set 
\begin{equation*}
A_{n,k}=\int_{S_{n,k}}J_{f}\,dtd\alpha .
\end{equation*}
The lemma implies that for every $n\geq 1$ the r.v. $A_{n,0}$, \dots , $%
A_{n,2^{n+1}-1}$ have the same law, which is also equal to the law of $\frac{%
1}{2}A_{n-1,0}$, \dots , $1/2A_{n,2^{n}-1}$. Therefore, setting $%
a_{n}=E\left[ A_{n,k}\right] $ (it may be infinite and does not depend on $k$%
), we have $a_{n}=a_{n-1}/2$ and (the integrand is positive so we
may interchange the operations) 
\begin{eqnarray*}
E\int_{S_{n}}J_{f}\,dtd\alpha
&=&E\sum_{k=0}^{2^{n+1}-1}A_{n,k}=\sum_{k=0}^{2^{n+1}-1}a_{n} \\
&=&2^{n+1}a_{n}=2^{n}a_{n-1}=...=2a_{0}\text{.}
\end{eqnarray*}
The claim of the lemma follows at once.
\end{proof}

\begin{acknowledgement}
The first named author (F.F.) thanks F. Russo for his insistence on
Lyons-Zheng processes in various discussions and 
C. Grisanti for very useful bibliographical advices.
\end{acknowledgement}

\end{document}